\documentclass{amsart}%
\usepackage{amsmath}
\usepackage{amscd}
\usepackage{amsfonts}
\usepackage{amssymb}
\usepackage{graphicx}%
\newtheorem{theorem}{Theorem}
\theoremstyle{plain}

\newtheorem{corollary}{Corollary}

\newtheorem{definition}{Definition}
\newtheorem{example}{Example}

\newtheorem{lemma}{Lemma}

\newtheorem{proposition}{Proposition}

\numberwithin{equation}{section}

\begin{document}
\title{A Diagonal on the Associahedra}
\author{Samson Saneblidze}
\address{A. Razmadze Mathematical Institute\\
Georgian Academy of Sciences\\
M. Aleksidze st., 1\\
0193 Tbilisi, Georgia}
\email{SANE@rmi.acnet.ge <mailto:SANE@rmi.acnet.ge> }
\author{Ronald Umble$^{1}$}
\address{Department of Mathematics\\
Millersville University of Pennsylvania\\
Millersville, PA. 17551}
\email{ron.umble@millersville.edu
%TCIMACRO{\TEXTsymbol{<}}%
%BeginExpansion
$<$%
%EndExpansion
mailto:ron.umble@millersville.edu%
%TCIMACRO{\TEXTsymbol{>} }%
%BeginExpansion
$>$
%EndExpansion
}
\thanks{$^{1}$This research funded in part by a Millersville University faculty
research grant.}

\begin{abstract}
Let $C_{\ast}\left(  K\right)  $ denote the cellular chains on the Stasheff
associahedra. We construct an explicit combinatorial diagonal $\Delta:C_{\ast
}(K)\rightarrow C_{\ast}(K)\otimes C_{\ast}(K);$ consequently, we obtain an
explicit diagonal on the ${A}_{\infty}$-operad. We apply the diagonal $\Delta$
to define the tensor product of $A_{\infty}$-(co)algebras in maximal generality.

\end{abstract}
\date{November 8, 2000; revised June 5, 2004}
\maketitle

\section{Introduction}

Let $C_{\ast}\left(  K\right)  $ denote the cellular chains on the disjoint
union of the Stasheff associahedra $\left\{  K_{n}\right\}  _{n\geq2}.$ In
this paper we construct an explicit combinatorial diagonal $\Delta:C_{\ast
}(K)\rightarrow C_{\ast}(K)\otimes C_{\ast}(K)$ based on a direct
decomposition of the top dimensional cells of $K.$ This leads to an explicit
diagonal on the ${A}_{\infty}$-operad and solves a long-standing problem. We
apply the diagonal $\Delta$ to define the tensor product of $A_{\infty}%
$-(co)algebras in maximal generality. We also include an appendix in which we
define an associahedral set $\mathcal{K}$ and lift $\Delta$ to a diagonal on
the chain complex of $\mathcal{K}$.

We mention that Chapoton \cite{Chapoton1}, \cite{Chapoton2} constructed a
diagonal on $C_{\ast}\left(  K\right)  $ of the form $\Delta:C_{\ast}\left(
K_{n}\right)  \rightarrow\bigoplus_{i+j=n}C_{\ast}\left(  K_{i}\right)
\otimes C_{\ast}\left(  K_{j}\right)  ,$ which coincides with the diagonal of
Loday and Ronco \cite{Loday} in dimension zero. Whereas Chapoton's diagonal is
primitive on generators, our diagonal is defined by geometrically decomposing
the generators. Thus the two are totally different.

\section{The Stasheff Associahedra}

In his seminal papers of 1963, J. Stasheff \cite{Stasheff} constructs the
associahedra\linebreak$\left\{  K_{n+2}\right\}  _{n\geq0}$ as follows: Let
$K_{2}=\ast;$ if $K_{n+1}$ has been constructed, let%
\[
L_{n+2}=\bigcup\limits_{\substack{r+s=n+3\\1\leq k\leq n-s+3}}\left(
K_{r}\times K_{s}\right)  _{k}%
\]
and define $K_{n+2}=CL_{n+2},$ i.e., the cone on $L_{n+2}.$ The associahedron
$K_{n+2}$ is an $n$-dimensional polyhedron, which serves as a parameter space
for homotopy associativity in $n+2$ variables. The top dimensional face of
$K_{n+2}$ corresponds to a pair of level 1 parentheses enclosing all $n+2$
indeterminants; each component $(K_{n-\ell+2}\times K_{\ell+1})_{i+1}$ of
$\partial K_{n+2}$ corresponds to a pair of level 2 parentheses enclosing
$\ell+1$ indeterminants beginning with the $\left(  i+1\right)  ^{st}$. We
denote this parenthesization by
\[
d_{\left(  i,\ell\right)  }=\left(  x_{1}\cdots\left(  x_{i+1}\cdots
x_{i+\ell+1}\right)  \cdots x_{n+2}\right)
\]
and refer to the inner and outer parentheses as the \textit{first }and
\textit{last }pair, respectively. Note that indices $i$ and $\ell$ are
constrained by%
\[
\left\{
\begin{array}
[c]{ll}%
0\leq i\leq n & \\
1\leq\ell\leq n, & i=0\\
1\leq\ell\leq n+1-i, & 1\leq i\leq n
\end{array}
\right\}  .
\]
Thus, there is a one-to-one correspondence between $\left(  n-1\right)
$-faces of $K_{n+2}$ and parenthesizations $d_{\left(  i,\ell\right)  }$ of
$n+2$ indeterminants.

Alternatively, $K_{n+2}$ can be realized as a subdivision of the standard
$n$-cube $I^{n}$ in the following way: Let $\epsilon=0,1.$ Label the endpoints
of $K_{3}=\left[  0,1\right]  $ via $\epsilon\leftrightarrow d_{\left(
\epsilon,1\right)  }.$ For $1\leq i\leq n,$ let $e_{i,\epsilon}^{n-1}$ denote
the $\left(  n-1\right)  $-face $(x_{1},\ldots,x_{i-1},\epsilon,x_{i+1}%
,\ldots,x_{n})\subset I^{n}$ and obtain $K_{4}$ from $K_{3}\times I=I^{2}$ by
subdividing the edge $e_{1,1}^{1}$ as the union of intervals $1\times
I_{0,1}\cup1\times I_{1,\infty}.$ Label the edges of $K_{4}\ $as follows:
$e_{i,0}^{1}\leftrightarrow d_{\left(  0,i\right)  };$ $e_{2,1}^{1}%
\leftrightarrow d_{\left(  2,1\right)  };$ $1\times I_{0,1}\leftrightarrow
d_{\left(  1,1\right)  };$ and $1\times I_{1,\infty}\leftrightarrow d_{\left(
1,2\right)  }$ (see Figure 1). Now for $0\leq i\leq j\leq\infty,$ let
$I_{i,j}$ denote the subinterval $\left[  \left(  2^{i}-1\right)
/2^{i},\left(  2^{j}-1\right)  /2^{j}\right]  \subset I,$ where $\left(
2^{\infty}-1\right)  /2^{\infty}$ is defined to be $1.$ For $n>2,$ assume that
$K_{n+1}$ has been constructed and obtain $K_{n+2}$ from $K_{n+1}\times
I\approx I^{n}$ by subdividing the $\left(  n-1\right)  $-faces $d_{\left(
i,n-i\right)  }\times I\ $as unions $d_{\left(  i,n-i\right)  }\times
I_{0,i}\cup d_{\left(  i,n-i\right)  }\times I_{i,\infty},$ $0<i<n$. Label the
$\left(  n-1\right)  $-faces of $K_{n+2}$ as follows:\vspace{0.2in}

$\hspace{0.2in}\hspace{0.2in}\hspace{0.2in}%
\begin{tabular}
[c]{l|ll}%
\textbf{Face of} $K_{n+2}$ & \textbf{Label} & \\\hline
&  & \\
$e_{\ell,0}^{n-1}$ & $d_{\left(  0,\ell\right)  },$ & $1\leq\ell\leq n$\\
&  & \\
$e_{n,1}^{n-1}$ & $d_{\left(  n,1\right)  ,}$ & \\
&  & \\
$d_{\left(  i,\ell\right)  }\times I$ & $d_{\left(  i,\ell\right)  ,}$ &
$1\leq\ell<n-i,\text{ }0<i<n-1$\\
&  & \\
$d_{\left(  i,n-i\right)  }\times I_{0,i}$ & $d_{\left(  i,n-i\right)  ,}$ &
$0<i<n$\\
&  & \\
$d_{\left(  i,n-i\right)  }\times I_{i,\infty}$ & $d_{\left(  i,n-i+1\right)
},$ & $0<i<n$%
\end{tabular}
$\vspace{0.2in}

\noindent In Figure 2 we have labeled the $2$-faces of $K_{5}$ that are
visible from the viewpoint of the diagram.

Compositions $d_{(i_{m},\ell_{m})}\dotsm d_{(i_{2},\ell_{2})}d_{(i_{1}%
,\ell_{1})}$ denote a successive insertion of $m+1$ pairs of parentheses into
$n+2$ indeterminants as follows: Given $d_{(i_{r},\ell_{r})}\dotsm
d_{(i_{1},\ell_{1})},$ $1\leq r<m,$ regard each pair of level 2 parentheses
and its contents as a single indeterminant and apply $d_{\left(  i_{r+1}%
,\ell_{r+1}\right)  }.$ Conclude by inserting a last pair enclosing
everything. Note that each parenthesization can be expressed as a unique
composition $d_{(i_{m},\ell_{m})}\dotsm d_{(i_{1},\ell_{1})}$ with
$i_{r+1}\leq i_{r}$ for $1\leq r<m,$ in which case the parentheses inserted by
$d_{\left(  i_{r+1},\ell_{r+1}\right)  }$ begin at or to the left of the pair
inserted by $d_{\left(  i_{r},\ell_{r}\right)  }$. Such compositions are said
to have \textit{first fundamental form}. Thus for $0\leq k<n,$ the $k$-faces
of $K_{n+2}$ lie in one-to-one correspondence with compositions $d_{(i_{n-k}%
,\ell_{n-k})}\dotsm d_{(i_{1},\ell_{1})}$ in first fundamental form. The two
extremes with $m$ pairs of parentheses inserted as far to the left and right
as possible, are respectively denoted by%
\[
d_{(0,\ell_{m})}\dotsm d_{(0,\ell_{1})}\hspace{0.1in}\text{and}\hspace
{0.1in}d_{(i_{m},i_{m-1}-i_{m})}\dotsm d_{(i_{1},i_{0}-i_{1})},
\]
where $i_{0}=n+1$ and $i_{r+1}<i_{r},\ 0\leq r<m.$ In particular, the $n$-fold
compositions%
\[
d_{(0,1)}\dotsm d_{(0,1)}\text{\hspace{0.1in}and\hspace{0.1in}}d_{(1,1)}\dotsm
d_{(n,1)}%
\]
denote the extreme full parenthesizations\ of $n+2$ indeterminants. When
$m=0,$ define $d_{(i_{m},\ell_{m})}\dotsm d_{(i_{1},\ell_{1})}=Id$.

\hspace*{0.4in}\setlength{\unitlength}{0.00043in}\begin{picture}
(200,1200)
\thicklines\put(3401,239){\line( 1,
0){1800}} \put(3401,239){\makebox(0,0){$\bullet$}}
\put(5201,239){\line( 0, 1){ 0}}
\put(5201,239){\makebox(0,0){$\bullet$}} \put(5201,239){\line(
0,-1){1800}} \put(5201,-661){\makebox(0,0){$\bullet$}}
\put(5201,-1561){\line(-1, 0){1800}}
\put(5201,-1561){\makebox(0,0){$\bullet$}} \put(3401,-1561){\line(
0, 1){1800}} \put(3401,-1561){\makebox(0,0){$\bullet$}}
\put(3401,239){\line( 0, 1){ 0}}
\put(4300,-680){\makebox(0,0){$K_{4}$}} \put(3110,-1861){\makebox
(0,0){$(0,0)$}} \put(3110,464){\makebox(0,0){$(0,1)$}}
\put(5500,464){\makebox (0,0){$(1,1)$}}
\put(5500,-1861){\makebox(0,0){$(1,0)$}} \put
(2910,-680){\makebox(0,0){$d_{(0,1)}$}} \put(4350,530){\makebox
(0,0){$d_{(2,1)}$}} \put(5700,-211){\makebox(0,0){$d_{(1,2)}$}}
\put (5700,-1060){\makebox(0,0){$d_{(1,1)}$}}
\put(4350,-1890){\makebox (0,0){$d_{(0,2)}$}}
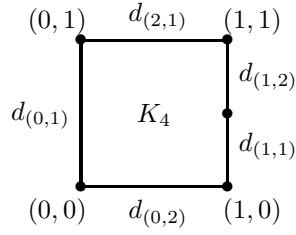
\end{picture}\vspace{1in}\newline\hspace*{1.2in}Figure 1: $K_{4}$ as a
subdivision of $K_{3}\times I.\ $

\hspace*{1.5in}\setlength{\unitlength}{0.00025in}\begin{picture}
(200,1800)
\thicklines\put(2601,-2761){\line(
0,-1){4800}} \put(2601,-2761){\makebox(0,0){$\bullet$}}
\put(2601,-7561){\line( 1, 0){4800}}
\put(2601,-7561){\makebox(0,0){$\bullet$}} \put(7401,-7561){\line(
0, 1){4800}} \put(7401,-7561){\makebox(0,0){$\bullet$}}
\put(7401,-2761){\line (-1, 0){4800}}
\put(7401,-2761){\makebox(0,0){$\bullet$}} \put
(2601,-3961){\line( 1, 0){4800}}
\put(2601,-3961){\makebox(0,0){$\bullet$}}
\put(2601,-2761){\line(-5, 4){3000}}
\put(2601,-2761){\makebox(0,0){$\bullet$ }} \put(-399,-361){\line(
0,-1){4800}} \put(-399,-361){\makebox(0,0){$\bullet $}}
\put(-399,-5161){\line( 5,-4){3000}} \put(-399,-5161){\makebox
(0,0){$\bullet$}} \put(-399,-361){\line( 1, 0){4800}}
\put(-399,-361){\makebox (0,0){$\bullet$}} \put(4401,-361){\line(
5,-4){3000}} \put(4401,-361){\makebox (0,0){$\bullet$}}
\put(1101,-1561){\line( 0,-1){4800}} \put
(1101,-1561){\makebox(0,0){$\bullet$}} \put(2601,-5161){\line(-5,
4){1500}} \put(2601,-5161){\makebox(0,0){$\bullet$}}
\put(-399,-5161){\line( 1, 0){1200}} \put(1401,-5161){\line( 1,
0){900}} \put(2901,-5161){\line( 1, 0){1500}}
\put(4401,-5161){\line( 0, 1){900}} \put(4401,-5161){\makebox
(0,0){$\bullet$}} \put(4401,-2461){\line( 0, 1){2100}} \put
(1100,-3970){\makebox(0,0){$\bullet$}} \put(7401,-3961){\makebox
(0,0){$\bullet$}} \put(4401,-5161){\line(5,-4){3000}} \put
(1100,-6350){\makebox(0,0){$\bullet$}} \put(-1350,-5500){\makebox
(0,0){$(1,0,0)$}} \put(4500,200){\makebox(0,0){$(0,0,1)$}} \put
(8150,-8161){\makebox(0,0){$(0,1,0)$}} \put(3226,-1586){\makebox
(0,0){$d_{(3,1)}$}} \put(326,-3061){\makebox(0,0){$d_{(1,1)}$}}
\put (1900,-3361){\makebox(0,0){$d_{(1,3)}$}}
\put(1900,-5711){\makebox (0,0){$d_{(1,2)}$}}
\put(4300,-6200){\makebox(0,0){$d_{(2,1)}$}}
\put(4700,-3286){\makebox(0,0){$d_{(2,2)}$}}
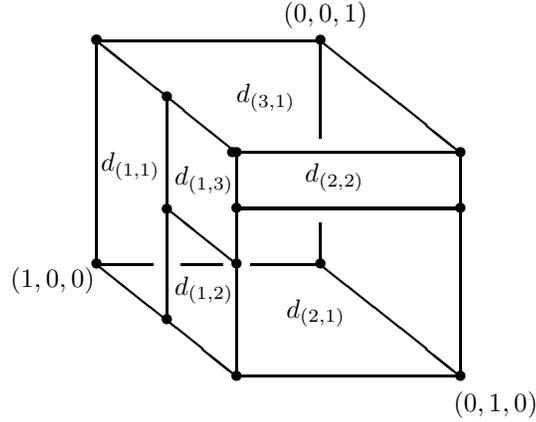
\end{picture}\vspace{2.2in}\newline\hspace*{1.2in}Figure 2: $K_{5}$ as a
subdivision of $K_{4}\times I.\vspace*{0.2in}$

Alternatively, each face of $K_{n+2}$ can be represented as a planar rooted
tree (PRT)\ with $n+2$ leaves; its leaves correspond to indeterminants and its
nodes correspond to pairs of parentheses. Let $T_{n+2}$ denote the PRT with
$n+2$ leaves attached to the root at a single node $N_{0},$ called the
\textit{root node} (see Figure 3). The leaves correspond to a single pair of
parentheses enclosing all $n+2$ indeterminants. Now given an arbitrary PRT
$T$, consider a node $N$ of valence $r+1\geq4$. Choose a neighborhood $U$ of
$N$ that excludes the other nodes of $T$ and note that $T_{r}\subseteq U\cap
T.$ Labeling from left to right, index the leaves of $T_{r}$ from $1$ to $r$
as in Figure 3. Perform an $\left(  i,\ell\right)  $-\textit{surgery at node
}$N$ in the following way: Remove leaves $i+1,\ldots,i+\ell+1$ of $T_{r},$
reattach them at a new node $N^{\prime}\neq N$ and graft in a new branch
connecting $N$ to $N^{\prime}$ (see Figure 4). Now let $n\geq1$. Given a
parenthesization $d_{\left(  i,\ell\right)  }$ of $n+2$ indeterminants, obtain
the PRT $T_{n+2}^{\left(  i,\ell\right)  }$ from $T_{n+2}$ by performing an
$\left(  i,\ell\right)  $-surgery at the root node $N_{0}$ as shown in Figure
4. Inductively, given a parenthesization $d_{(i_{m},\ell_{m})}\dotsm
d_{(i_{1},\ell_{1})}$ of $n+2$ indeterminants expressed as a composition in
first fundamental form, construct the corresponding PRT $T_{n+2}^{\left(
i_{1},\ell_{1}\right)  ,\ldots,\left(  i_{m},\ell_{m}\right)  }$ as follows:
Assume that $T_{n+2}^{\left(  i_{1},\ell_{1}\right)  ,\ldots,\left(
i_{r},\ell_{r}\right)  }$ with nodes $N_{0},\ldots,N_{r}$ has been constructed
for some $1\leq r<m$ and note that the root node $N_{0}$ has valence
$n+3-\ell_{1}-\cdots-\ell_{r}.$ Perform an $(i_{r+1},\ell_{r+1})$-surgery
at\textit{\ }$N_{0}$ and obtain $T_{n+2}^{\left(  i_{1},\ell_{1}\right)
,\ldots,\left(  i_{r+1},\ell_{r+1}\right)  }$ containing a new node $N_{r+1}$
and a new branch connecting $N_{0}$ to $N_{r+1}.$ Finally, define
$T_{n+2}^{\left(  i_{1},\ell_{1}\right)  ,\ldots,\left(  i_{m},\ell
_{m}\right)  }=T_{n+2}$ when $m=0$ and obtain a one-to-one correspondence
between $k$-faces of $K_{n+2},$ $0\leq k\leq n$ and PRT's $T_{n+2}^{\left(
i_{1},\ell_{1}\right)  ,\ldots,\left(  i_{n-k},\ell_{n-k}\right)  }$
consisting of $n-k+1$ nodes and $n+2$ leaves. In particular, each vertex of
$K_{n+2}$ corresponds to a planar binary rooted tree $T_{n+2}^{\left(
i_{1},\ell_{1}\right)  ,\ldots,\left(  i_{n},\ell_{n}\right)  }$ (see Figure 5).

\begin{center}
\setlength{\unitlength}{0.00015in} \begin{picture}
(2775,3685)(2926,-1038)
\thicklines\put(4350,-2580){\makebox(0,0){$\bullet$ }}
%node N0
\put(4200,-2700){\line(-1, 1){2700}}
%leaf1
\put(4200,-2700){\line(-1, 2){1350}}
%leaf2
\put(4200,0){\line(0, -1){4100}}
%leaf3 & root
\put(5500,-500){\makebox(0,0){$\dots$ }}
\put(4200,-2700){\line(1, 1){2700}}
%leaf5
\put(3500,-3250){\makebox(0,0){$N_{0}$ }}  \put(4350,900){\makebox
(0,0){$1\hspace*{0.2in} 2\hspace*{0.2in} 3\hspace*{0.05in}  \cdots
\hspace*{0.05in}n+2$}}
\end{picture}\vspace{0.5in}

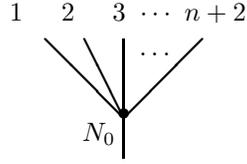
Figure 3: The corolla $T_{n+2}$\vspace{0.5in}\linebreak%
\setlength{\unitlength}{0.0003in}\begin{picture}
(2775,3685)(2926,-1038)
\thicklines\put(4285,-2650){\makebox(0,0){$\bullet$ }}
%node N0
\put(3700,-3100){\makebox(0,0){$N_{0}$ }}
\put(4285,1050){\makebox (0,0){$\bullet$ }}
%node N'
\put(3700,500){\makebox(0,0){$N^{\prime}$ }}
\put(4200,-2700){\line(-1, 1){2700}}
%leaf1
\put(4200,-2700){\line(-1, 2){1350}}
%leaf2
\put(4200,1000){\line(0, -1){5000}}
%leaf3 & root
\put(4200,-2700){\line(1, 2){1350}} \put(4200,-2700){\line(1,
1){2700}}
%leaf5
\put(4850,900){\makebox(0,0){$1\hspace*{0.1in}
\cdots\hspace*{0.15in}i \hspace*{0.7in}  i+\ell+2 \cdots n+2$}}
\put(4300,3300){\makebox(0,0){$i+1 \hspace*{0.2in}
\cdots\hspace*{0.2in} i+\ell+1 $}}  \put(4200,1000){\line(-1,
1){1500}}
%leaf1
\put(4200,1000){\line(-1, 2){750}}
%leaf2
\put(4300,2300){\makebox(0,0){$\dots$ }}
\put(2500,-300){\makebox (0,0){$\dots$ }}
\put(6100,-300){\makebox(0,0){$\dots$ }} \put (4200,1000){\line(1,
2){750}} \put(4200,1000){\line(1, 1){1500}}
%leaf5
\end{picture}\vspace{1in}

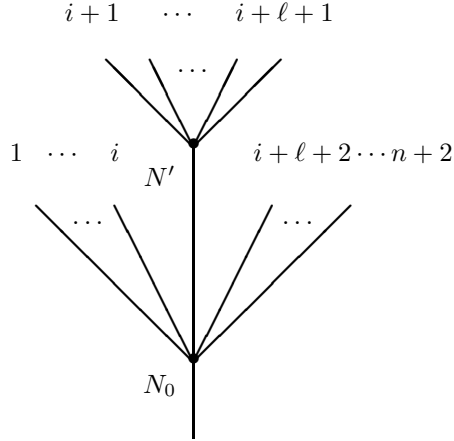
Figure 4: The PRT $T_{n+2}^{\left(  i,\ell\right)  }$\vspace{0.2in}

\setlength{\unitlength}{0.0002in}\begin{picture}
(2775,3685)(2926,-1038) \thicklines\put(-1000,-561){\line(0,
-1){900} }
%root
\put(-1000,-561){\makebox(0,0){$\bullet$}}
%node N0
\put(-1000,-561){\line(-1, 1){2700}}
%left branch
\put(-1000,-561){\line(1, 1){2700}}
%right branch
\put(-1940,339){\line(1, 1){1800}}
%leaf 3
\put(-2840,1239){\line(1, 1){900}}
%leaf 2
\put(-1900,339){\makebox(0,0){$\bullet$}}
%node N1
\put(-2800,1239){\makebox(0,0){$\bullet$}}
%node N2
\put(-1000,-2800){\makebox(0,0){$T_{4}^{(0,1),(0,1)}
\leftrightarrow$}}
\put(-1200,-4000){\makebox(0,0){$d_{(0,1)}d_{(0,1)} = (((
\bullet\bullet) \bullet) \bullet)$}}   \put(9500,-561){\line(0,
-1){900}}
%root
\put(9500,-561){\makebox(0,0){$\bullet$}}
%node N0
\put(9500,-561){\line(-1, 1){2700}}
%left branch
\put(9500,-561){\line(1, 1){2700}}
%right branch
\put(10360,339){\line(-1, 1){1800}}
%leaf 2
\put(11260,1239){\line(-1, 1){900}}
%leaf 3
\put(10400,339){\makebox(0,0){$\bullet$}}
%node N1
\put(11300,1239){\makebox(0,0){$\bullet$}}
%node N2
\put(9500,-2800){\makebox(0,0){$T_{4}^{(2,1),(1,1)}
\leftrightarrow$}}
\put(9300,-4000){\makebox(0,0){$d_{(1,1)}d_{(2,1)} = ( \bullet(
\bullet( \bullet\bullet)))$}}
\end{picture}
\vspace{0.8in}

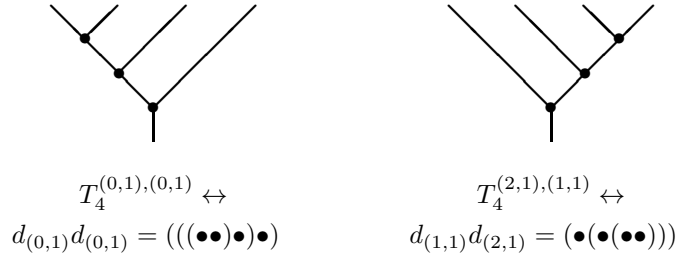
Figure 5: Some binary PRT's.\vspace{0.1in}
\end{center}

Now given a $k$-face $a_{k}\subseteq K_{n+2},$ $k>0,$ consider the two
vertices of $a_{k}$ at which parentheses are shifted as far to the left and
right as possible; we refer to these vertices as the minimal and maximal
vertices of $a_{k},$ and denote them by $a_{k}^{\min}$ and $a_{k}^{\max},$
respectively. In particular, the minimal and maximal vertices of $K_{n+2}$ are
the origin and the vertex of $I^{n}$ diagonally opposite to it, i.e.,
\[
K_{n+2}^{\min}\leftrightarrow\left(  0,0,\ldots,0\right)  \text{ and }%
K_{n+2}^{\max}\leftrightarrow\left(  1,1,\ldots,1\right)  ;
\]
the respective binary trees in Figure 5 correspond to $K_{4}^{\min}$ and
$K_{4}^{\max}$. Given a representation $T_{n+2}^{\left(  i_{1},\ell
_{1}\right)  ,\ldots,\left(  i_{n-k},\ell_{n-k}\right)  }$ of $a_{k}$,
construct the minimal (resp., maximal) tree of $a_{k}$ by replacing each node
of valence $r\geq4$ with the planar binary rooted tree representing
$K_{r-1}^{\min}$ (resp., $K_{r-1}^{\max}$). Note that $a_{k}^{\min}$ and
$a_{k}^{\max}$ determine $a_{k}$ since their convex hull is a diagonal of
$a_{k}$. But we can say more.

When a composition of face operators $d_{(i_{m},\ell_{m})}\dotsm
d_{(i_{1},\ell_{1})}$ is defined we refer to the sequence of lower indices
$I=\left(  i_{1},\ell_{1}\right)  ,\ldots,$ $\left(  i_{m},\ell_{m}\right)  $
as an \textit{admissible sequence of length} $m;$ if $d_{(i_{m},\ell_{m}%
)}\dotsm d_{(i_{1},\ell_{1})}$ has first fundamental form we refer to the the
sequence $I$ as a \textit{type I sequence of length} $m$. The set of all
planar binary rooted trees
\[
Y_{n+2}=\left\{  T_{n+2}^{I}\text{ }\left\vert \text{ }I\text{ is a type I
sequence of length }n\right.  \right\}
\]
is a poset with partial ordering defined as follows (cf. \cite{Pallo}): Say
that $T_{n+2}^{I_{p}}\leq T_{n+2}^{I_{q}}$ if there is an edge-path in
$K_{n+2}$ from vertex $T_{n+2}^{I_{p}}$ to vertex $T_{n+2}^{I_{q}}$ along
which parentheses shift strictly to the right. This partial ordering can be
expressed geometrically in terms of the following operation on tress: Let
$N_{0}$ denote the root node of $T_{n+2}^{I}$ and let $N$ be a node joining
some left branch $L$ and a right branch or leaf $R$ in $T_{n+2}^{I}.$ Let
$N_{L}$ denote the node on $L$ immediately above $N$. A \textit{right-shift
through node }$N$ repositions $N_{L}$ either at the midpoint of leaf $R$ or
midway between $N$ and the node immediately above it. Then $T_{n+2}^{I_{p}%
}\leq T_{n+2}^{I_{q}}$ if there is a \textit{right-shift sequence} of planar
binary rooted trees $\left\{  T_{n+2}^{I_{r}}\right\}  _{p\leq r\leq q},$
i.e., for each $r<q,$ tree $T_{n+2}^{I_{r+1}}$ is obtained from $T_{n+2}%
^{I_{r}}$ by a right-shift through some node in $T_{n+2}^{I_{r}}$ (see Figure
6).\vspace{0.2in}

\begin{center}
\setlength{\unitlength}{0.0002in} \begin{picture}
(2775,3685)(2926,-1038) \thicklines\put(-3500,-561){\line(0,
-1){900} }
%root
\put(-3500,-561){\makebox(0,0){$\bullet$}}
%node N0
\put(-3500,-561){\line(-1, 1){2700}}
%left branch
\put(-3500,-561){\line(1, 1){2700}}
%right branch
\put(-4440,339){\line(1, 1){1800}}
%leaf 3
\put(-5340,1239){\line(1, 1){900}}
%leaf 2
\put(-4400,339){\makebox(0,0){$\bullet$}}
%node N1
\put(-5300,1239){\makebox(0,0){$\bullet$}}
%node N2
\put(500,350){\makebox(0,0){$\leq$}}  \put(4250,-561){\line(0,
-1){900}}
%root
\put(4250,-561){\makebox(0,0){$\bullet$}}
%node N0
\put(4250,-561){\line(-1, 1){2700}}
%left branch
\put(4250,-561){\line(1, 1){2700}}
%right branch
\put(2410,1239){\line(1, 1){900}}
%leaf 2
\put(6010,1239){\line(-1, 1){900}}
%leaf 3
\put(2450,1239){\makebox(0,0){$\bullet$}}
%node N1
\put(6050,1239){\makebox(0,0){$\bullet$}}
%node N2
\put(8225,350){\makebox(0,0){$\leq$}}  \put(12000,-561){\line(0,
-1){900}}
%root
\put(12000,-561){\makebox(0,0){$\bullet$}}
%node N0
\put(12000,-561){\line(-1, 1){2700}}
%left branch
\put(12000,-561){\line(1, 1){2700}}
%right branch
\put(12860,339){\line(-1, 1){1800}}
%leaf 2
\put(13860,1239){\line(-1, 1){900}}
%leaf 3
\put(12900,339){\makebox(0,0){$\bullet$}}
%node N1
\put(13800,1239){\makebox(0,0){$\bullet$}}
%node N2
\end{picture}
\vspace{0.3in}

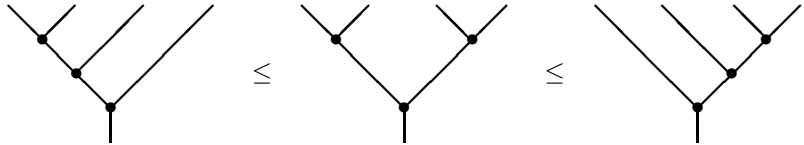
Figure 6: A right-shift sequence of planar binary trees\vspace{0.2in}
\end{center}

Let $I^{\prime}=(i_{1}^{\prime},\ell_{1}^{\prime}),...,(i_{m}^{\prime}%
,\ell_{m}^{\prime})$ be an admissible sequence of length $m>0$ and consider a
node $N^{\prime}$ distinct from the root node $N_{0}$ in the PRT
$T_{n+2}^{I^{\prime}}.$ Let $N$ denote the node immediately below $N^{\prime}$
and let $NN^{\prime}$ denote the branch from $N$ to $N^{\prime};$ we refer to
the quotient space $T_{n+2}^{I}=T_{n+2}^{I^{\prime}}/NN^{\prime}$\ as the
$\left(  N,N^{\prime}\right)  $\textit{-contraction} \textit{of}
$T_{n+2}^{I^{\prime}}.$ Now given a type I sequence $J^{\prime}$ of length
$n$, consider the planar binary tree $T_{n+2}^{J^{\prime}}$ and let
$T_{n+2}^{J}$ be the PRT obtained from $T_{n+2}^{J^{\prime}}$ by some sequence
of $k$ successive $\left(  N,N^{\prime}\right)  $-contractions. The subposet
$Y_{n+2}^{J}\subseteq Y_{n+2}$ of all planar binary rooted trees from which
$T_{n+2}^{J}$ can be so obtained is exactly the poset of vertices of the
$k$-face $a_{k}\subseteq K_{n+2}$ represented by $T_{n+2}^{J}.$ In this way,
we may regard $a_{k}$ as the geometric realization of $Y_{n+2}^{J}$ just as we
regard a $k$-face of the standard $n$-simplex as the geometric realization of
a $\left(  k+1\right)  $-subset of a linearly ordered $\left(  n+1\right)
$-set. In particular, $K_{n+2}$ is the geometric realization of $Y_{n+2}$.

We summarize the discussion above as a proposition:

\begin{proposition}
For $0\leq k\leq n,$ the following correspondences preserve combinatorial
structure:%
\[%
\begin{array}
[c]{lll}%
\left\{  k\text{\textit{-faces of} }K_{n+2}\right\}  & \leftrightarrow &
\left\{
\begin{array}
[c]{c}%
\left(  n-k\right)  \text{-fold c\textit{ompositions} of face}\\
\text{operators \textit{in first fundamental form}}%
\end{array}
\right\}  \vspace{0.1in}\\
& \leftrightarrow & \left\{
\begin{array}
[c]{c}%
\text{\textit{Planar rooted trees with}}\\
n-k+1\ \text{\textit{nodes and }}n+2\ \text{\textit{leaves}}%
\end{array}
\right\}  \vspace{0.1in}\\
& \leftrightarrow & \left\{
\begin{array}
[c]{c}%
\text{S\textit{ubposets of planar binary rooted trees }}Y_{n+2}^{J}\\
\text{where }J\text{ is a type I sequence of length }n-k
\end{array}
\right\}  .
\end{array}
\]

\end{proposition}

\section{A Diagonal $\Delta$ on $C_{\ast}\left(  K_{n}\right)  $\label{diag}}

For notational simplicity, we suppress upper indices $q_{2},\ldots,q_{m}$ in a
composition $d_{(i_{m},\ell_{m})}^{q_{m}}\dotsm d_{(i_{2},\ell_{2})}^{q_{2}%
}d_{(i_{1},\ell_{1})}^{q_{1}}$ when $q_{j+1}=q_{j}+1$ for all $j\geq1;$ if, in
addition, $q_{1}=1,$ we suppress all $q_{i}$'s.

\begin{definition}
Let $m\geq2.$ A sequence of lower indices $I=\left(  i_{1},\ell_{1}\right)
,\ldots,$ $\left(  i_{m},\ell_{m}\right)  $ is \underline{admissible} whenever
the composition of face operators $d_{(i_{m},\ell_{m})}^{q_{m}}\dotsm
d_{(i_{1},\ell_{1})}^{q_{1}}$ is defined. The sequence $I$ is a type I (resp.
type II) sequence if $I$ is \underline{admissible} and $i_{k}\geq i_{k+1}$
(resp. $i_{k}\leq i_{k+1}+\ell_{k+1}$) for $1\leq k<m$. The empty sequence
$\left(  m=0\right)  $ and sequences of length $1$ $\left(  m=1\right)  $ are
sequences of types I and II. A composition of face operators $d_{(i_{m}%
,\ell_{m})}\dotsm d_{(i_{2},\ell_{2})}d_{(i_{1},\ell_{1})}^{s}$ has
\underline{first (resp. second) fundamental form} if $\left(  i_{1},\ell
_{1}\right)  ,\ldots,$ $\left(  i_{m},\ell_{m}\right)  $ is a type I (resp.
type II) sequence. When $m=0,$ the composition $d_{(i_{m},\ell_{m})}\dotsm
d_{(i_{2},\ell_{2})}d_{(i_{1},\ell_{1})}^{s}$ is defined to be the identity.
An element $b=d_{(i_{m},\ell_{m})}^{q_{m}}\dotsm d_{(i_{2},\ell_{2})}^{q_{2}%
}d_{(i_{1},\ell_{1})}^{q_{1}}(a)$ is expressed in \underline{first (resp.
second) fundamental form} as a face of $a$ if%
\[
d^{q_{m}}\dotsm d^{q_{2}}d^{q_{1}}=\left(  d\dotsm dd^{s_{m}}\right)
\dotsm\left(  d\dotsm dd^{s_{2}}\right)  \left(  d\dotsm dd^{s_{1}}\right)  ,
\]
where $s_{1}<s_{2}<\cdots<s_{m}$ and each composition $d_{(i_{s_{j}}%
,\ell_{s_{j}})}\dotsm d_{(i_{2},\ell_{2})}d_{(i_{1},\ell_{1})}^{s_{j}}$ has
first (resp. second) fundamental form.
\end{definition}

Face operators satisfy the following relations:%
\[%
\begin{array}
[c]{lll}%
d_{(i_{p},\ell_{p})}^{p}d_{(i_{q},\ell_{q})}^{q}=\vspace{0.1in}d_{(i_{q}%
,\ell_{q})}^{q+1}d_{(i_{p},\ell_{p})}^{p}, & p<q & (1)\\
d_{(i_{q+1},\ell_{q+1})}^{q+1}d_{(i_{q},\ell_{q})}^{q}=\vspace{0.1in}%
d_{(i_{q}-i_{q+1},\ell_{q})}^{q}d_{(i_{q+1},\ell_{q+1}+\ell_{q})}^{q}, &
i_{q+1}\leq i_{q}\leq i_{q+1}+\ell_{q+1} & (2)\\
d_{(i_{q+1},\ell_{q+1})}^{q+1}d_{(i_{q},\ell_{q})}^{q}=\vspace{0.1in}%
d_{(i_{q},\ell_{q})}^{q+1}d_{(i_{q+1}+\ell_{q},\ell_{q+1})}^{q}, &
i_{q}<i_{q+1}. & (3)
\end{array}
\]
\noindent Furthermore, every composition of face operators can be uniquely
transformed into first or second fundamental form by successive applications
of face relations (1) to (3). For example, when $n=2$, the following five face
operators relate $T_{4}\in K_{4}^{2}$ to the edges of the pentagon $K_{4}$:
\[%
\begin{array}
[c]{cccc}%
d_{\left(  0,2\right)  }\left(  T_{4}\right)  & \longmapsto & \left(  \left(
\bullet\bullet\bullet\right)  \bullet\right)  & \in K_{3}\times K_{2}\\
d_{\left(  1,2\right)  }\left(  T_{4}\right)  & \longmapsto & \left(
\bullet\left(  \bullet\bullet\bullet\right)  \right)  & \in K_{3}\times
K_{2}\\
d_{\left(  0,1\right)  }\left(  T_{4}\right)  & \longmapsto & \left(  \left(
\bullet\bullet\right)  \bullet\bullet\right)  & \in K_{2}\times K_{3}\\
d_{\left(  1,1\right)  }\left(  T_{4}\right)  & \longmapsto & \left(
\bullet\left(  \bullet\bullet\right)  \bullet\right)  & \in K_{2}\times
K_{3}\\
d_{\left(  2,1\right)  }\left(  T_{4}\right)  & \longmapsto & \left(
\bullet\bullet\left(  \bullet\bullet\right)  \right)  & \in K_{2}\times K_{3}.
\end{array}
\]
There are four compositions of face operators
\[
d_{\left(  i_{2},1\right)  }^{1}d_{\left(  i_{1},2\right)  }^{1}%
:K_{4}\rightarrow K_{3}\times K_{2}\rightarrow K_{2}\times K_{2}\times K_{2}%
\]
with $0\leq i_{1},i_{2}\leq1,$ and six compositions
\[
d_{\left(  i_{2},1\right)  }^{2}d_{\left(  i_{1},1\right)  }^{1}%
:K_{4}\rightarrow K_{2}\times K_{3}\rightarrow K_{2}\times K_{2}\times K_{2}%
\]
with $0\leq i_{1}\leq2$ and $0\leq i_{2}\leq1,$ which pair off via relations
(1) to (3) and relate $T_{4}$ to each of the five vertices of $K_{4}$:
\[%
\begin{array}
[c]{ccccc}%
d_{\left(  0,1\right)  }^{2}d_{\left(  0,1\right)  }^{1}(T_{4}) & = &
d_{\left(  0,1\right)  }^{1}d_{\left(  0,2\right)  }^{1}(T_{4}) & \longmapsto
& \left(  \left(  \left(  \bullet\bullet\right)  \bullet\right)
\bullet\right) \\
d_{\left(  0,1\right)  }^{2}d_{\left(  1,1\right)  }^{1}(T_{4}) & = &
d_{\left(  1,1\right)  }^{1}d_{\left(  0,2\right)  }^{1}(T_{4}) & \longmapsto
& \left(  \left(  \bullet\left(  \bullet\bullet\right)  \right)
\bullet\right) \\
d_{\left(  1,1\right)  }^{2}d_{\left(  1,1\right)  }^{1}(T_{4}) & = &
d_{\left(  0,1\right)  }^{1}d_{\left(  1,2\right)  }^{1}(T_{4}) & \longmapsto
& \left(  \bullet\left(  \left(  \bullet\bullet\right)  \bullet\right)
\right) \\
d_{\left(  1,1\right)  }^{2}d_{\left(  2,1\right)  }^{1}(T_{4}) & = &
d_{\left(  1,1\right)  }^{1}d_{\left(  1,2\right)  }^{1}(T_{4}) & \longmapsto
& \left(  \bullet\left(  \bullet\left(  \bullet\bullet\right)  \right)
\right) \\
d_{\left(  0,1\right)  }^{2}d_{\left(  2,1\right)  }^{1}(T_{4}) & = &
d_{\left(  1,1\right)  }^{2}d_{\left(  0,1\right)  }^{1}(T_{4}) & \longmapsto
& \left(  \left(  \bullet\bullet\right)  \left(  \bullet\bullet\right)
\right)  .
\end{array}
\]
These relations encode the fact that when inserting two pairs of parentheses
into a string of four variables either pair may be inserted first.

Given a type $I$ sequence $(j_{1},n_{1}+1),...,(j_{k},n_{k}+1)$ and a sequence
$(i_{1},...,i_{k+1})$ with $0\leq i_{q}\leq n_{q},\,1\leq q\leq k+1,$ define a
function of two variables
\begin{equation}
{j(q,r)}=\left\{
\begin{array}
[c]{ll}%
i_{q}+j_{q}+n_{q-1}+\cdots+n_{r}+q-r, & 1\leq r<q\\
i_{q}+j_{q}, & 1\leq r=q
\end{array}
\right.  \label{ksign1}%
\end{equation}
for $1\leq r\leq q\leq k+1$ (assuming $j_{k+1}=0$), and let%
\begin{equation}
\beta=\left\{
\begin{array}
[c]{ll}%
1, & q=1,\\
\max\limits_{1\leq r\leq q}\{r\,|\,j(q,r)\leq j_{r-1}\}, & q>1,
\end{array}
\right.  \label{ksign2}%
\end{equation}
(assuming $j_{0}=\infty.$)

\begin{definition}
\label{sign}For $k\geq0,$ let $I_{k}=(j_{1},n_{1}+1),(j_{2},n_{2}%
+1),\ldots,(j_{k},n_{k}+1)$ be a type I sequence. If $1\leq q\leq k+1,$ the
\underline{sign of the face} $b=d_{(i_{q},\ell_{q})}^{q}(T_{n+2}^{I_{k}})\in
C_{n-k-1}\left(  K_{n+2}\right)  $\ is defined to be $(-1)^{\epsilon
_{1}+\epsilon_{2}},$ where
\begin{equation}%
\begin{array}
[c]{l}%
\epsilon_{1}=\left(  i_{q}+1\right)  \ell_{q}+n_{1}+\cdots+n_{q-1},\\
\epsilon_{2}=\left\{
\begin{array}
[c]{ll}%
0, & 1\leq\beta=q,\\
(\ell_{q}-1)(n_{q-1}+\cdots+n_{\beta}), & 1\leq\beta<q,
\end{array}
\right.
\end{array}
\label{esign}%
\end{equation}
and $\beta$ is defined by (\ref{ksign2}).
\end{definition}

Let us construct an explicit diagonal on associahedra in terms of compositions
of face operators in first and second fundamental form; the formulas we obtain
also determine a DG coalgebra structure on $(C_{\ast}(K_{n+2}),d)$.

We begin with an overview of the geometric ideas involved. Let $0\leq q\leq n$
and let $I_{n-q}$ be a type I sequence. The $q$-dimensional generator
$a_{q}=T_{n+2}^{I_{n-q}}$ is associated with a face of $K_{n+2}$ corresponding
to $n+2$ indeterminants with $n-$ $q+1$ pairs of parentheses. Identify $a_{q}$
with its associated face of $K_{n+2}$ and consider the minimal and maximal
vertices $a_{q}^{\min}$ and $a_{q}^{\max}$of $a_{q}$. Define the primitive
terms of $\Delta T_{n+2}$ to be
\[
T_{n+2}^{\min}\otimes T_{n+2}+T_{n+2}\otimes T_{n+2}^{\max}.
\]
Let $0<p<p+q=n$ and consider distinct $p$-faces $b$ and $b^{\prime}$ of
$T_{n+2}.$ Say that $b\leq b^{\prime}$ if there is a path of $p$-faces from
$b$ to $b^{\prime}$ along which parentheses shift strictly to the right. Now
given a $q$-dimensional face $a_{q}$ of $T_{n+2}$ such that $a_{q}^{\min}\neq
T_{n+2}^{\min},$ there is a unique path of $p$-faces $b_{1}\leq b_{2}%
\leq\cdots\leq b_{r}$ with minimal length such that $T_{n+2}^{\min}%
=b_{1}^{\min}$ and $a_{q}^{\min}=b_{r}^{\max}.$ Up to sign, we define the
non-primitive terms of $\Delta T_{n+2}$ to be
\[
\sum\pm\text{ }b_{j}\otimes a_{q}.
\]

To visualize this, consider the edge $d_{\left(  1,2\right)  }\left(
T_{4}\right)  \in C_{1}\left(  K_{4}\right)  $ whose minimal vertex is the
point $\left(  1,\frac{1}{2}\right)  $ (see Figure 7)$.$ The edges $d_{\left(
0,2\right)  }\left(  T_{4}\right)  \leq d_{\left(  1,1\right)  }\left(
T_{4}\right)  $ form a path of minimal length from $\left(  0,0\right)  $ to
$\left(  1,\frac{1}{2}\right)  $. Consequently, $\Delta T_{4}$ contains the
non-primitive terms $\left\{  \left(  \pm d_{\left(  0,2\right)  }\pm
d_{\left(  1,1\right)  }\right)  \otimes d_{\left(  1,2\right)  }\right\}
\left(  T_{4}\otimes T_{4}\right)  .\vspace{0.1in}$

\begin{center}
\setlength{\unitlength}{0.0005in}\begin{picture}
(2775,2685)(2926,-2038) \thicklines\put(3401,239){\line( 1,
0){1800}} \put(3401,239){\makebox(0,0){$\bullet$}}
\put(3400,155){\makebox (0,0){$\blacktriangle$}}
\put(5200,155){\makebox(0,0){$\blacktriangle$}}
\put(5200,-745){\makebox(0,0){$\blacktriangle$}}
\put(5105,241){\makebox (0,0){$\blacktriangleright$}}
\put(5105,-1559){\makebox (0,0){$\blacktriangleright$}}
\put(5200,239){\line( 0, 1){ 0}} \put
(5200,239){\makebox(0,0){$\bullet$}} \put(5200,239){\line(
0,-1){1800}} \put(5200,-661){\makebox(0,0){$\bullet$}}
\put(5200,-1561){\line(-1, 0){1800}}
\put(5200,-1561){\makebox(0,0){$\bullet$}} \put(3400,-1561){\line(
0, 1){1800}} \put(3400,-1561){\makebox(0,0){$\bullet$}}
\put(3400,239){\line( 0, 1){ 0}}
\put(4300,-680){\makebox(0,0){$T_{4}$}} \put(3110,-1861){\makebox
(0,0){$T_{4}^{min}$}} \put(5500,464){\makebox(0,0){$T_{4}^{max}$}}
\put(2610,-680){\makebox(0,0){$d_{(0,1)}(T_{4})$}}
\put(4350,530){\makebox (0,0){$d_{(2,1)}(T_{4})$}}
\put(6000,-211){\makebox(0,0){$d_{(1,2)}(T_{4})$}}
\put(6000,-1060){\makebox(0,0){$d_{(1,1)}(T_{4})$}}
\put(4350,-1890){\makebox (0,0){$d_{(0,2)}(T_{4})$}}
\end{picture}$\vspace{0.1in}$

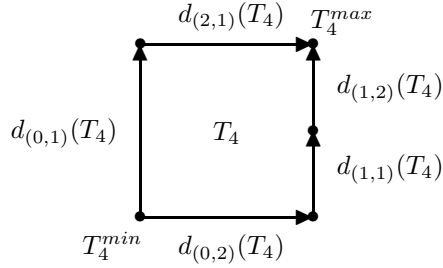
Figure 7: Edge paths from $T_{4}^{\min}$ to $T_{4}^{\max}$.$\vspace{0.1in}$
\end{center}

Precisely, for $T_{2}\in C_{\ast}\left(  K_{2}\right)  $ define $\Delta
T_{2}=T_{2}\otimes T_{2};$ inductively, assume that the map $\Delta:C_{\ast
}\left(  K_{i+2}\right)  \rightarrow C_{\ast}\left(  K_{i+2}\right)  \otimes
C_{\ast}\left(  K_{i+2}\right)  $ has been defined for all $i<n.$ For
$T_{n+2}\in C_{n}\left(  K_{n+2}\right)  $ define
\[
\Delta T_{n+2}=\sum_{0\leq p\leq p+q=n}\left(  -1\right)  ^{\epsilon}\text{
}d_{(i_{q}^{\prime},\ell_{q}^{\prime})}\cdots d_{(i_{1}^{\prime},\ell
_{1}^{\prime})}\left(  T_{n+2}\right)  \otimes d_{(i_{p},\ell_{p})}\cdots
d_{(i_{1},\ell_{1})}\left(  T_{n+2}\right)
\]
where
\[
\epsilon=\sum_{j=1}^{q}i_{j}^{\prime}(\ell_{j}^{\prime}+1)+\sum_{k=1}%
^{p}(i_{k}+k+p+1)\ell_{k},
\]
and lower indices $\left(  \left(  i_{1},\ell_{1}\right)  ,\ldots,\left(
i_{p},\ell_{p}\right)  ;\left(  i_{1}^{\prime},\ell_{1}^{\prime}\right)
,\ldots,\left(  i_{q}^{\prime},\ell_{q}^{\prime}\right)  \right)  $ range over
all solutions of the following system of inequalities:%

\begin{equation}
\hspace*{0.4in}\left\{
\begin{array}
[c]{ll}%
1\leq i_{j}<i_{j-1}\leq n+1\vspace{0.1in} & (1)\\
1\leq\ell_{j}\leq n+1-i_{j}-\ell_{\left(  j-1\right)  }\vspace{0.1in} & (2)\\
0\leq i_{k}^{\prime}\leq\min\limits_{o^{\prime}\left(  t_{k}\right)
<r<k}\left\{  i_{r}^{\prime},\text{ }i_{t_{k}}-\ell_{\left(  o^{\prime}\left(
t_{k}\right)  \right)  }^{\prime}\right\}  \vspace{0.1in} & (3)\\
1\leq\ell_{k}^{\prime}=\epsilon_{k}-i_{k}^{\prime}-\ell_{\left(  k-1\right)
}^{\prime} & (4)
\end{array}
\right\}  _{\substack{1\leq j\leq p\\1\leq k\leq q}}, \label{shuffle2}%
\end{equation}
where

$\hspace*{1.02in}\left\{  \epsilon_{1}<\cdots<\epsilon_{q}\right\}  =\left\{
1,\ldots,n\right\}  \setminus\left\{  i_{1},\ldots,i_{p}\right\}  ;$%
\vspace{0.1in}

$\hspace*{1.02in}\epsilon_{0}=\ell_{0}=\ell_{0}^{\prime}=i_{p+1}%
=i_{q+1}^{\prime}=0;$\vspace{0.1in}

$\hspace*{1.02in}i_{0}=i_{0}^{\prime}=\epsilon_{q+1}=\ell_{\left(  p+1\right)
}=\ell_{\left(  q+1\right)  }^{\prime}=n+1;$\vspace{0.1in}

$\hspace*{1.02in}\ell_{\left(  u\right)  }=%
%TCIMACRO{\tsum \nolimits_{j=0}^{u}}%
%BeginExpansion
{\textstyle\sum\nolimits_{j=0}^{u}}
%EndExpansion
\ell_{j}$ for $0\leq u\leq p+1;$\vspace{0.1in}

$\hspace*{1.02in}\ell_{\left(  u\right)  }^{\prime}=%
%TCIMACRO{\tsum \nolimits_{k=0}^{u}}%
%BeginExpansion
{\textstyle\sum\nolimits_{k=0}^{u}}
%EndExpansion
\ell_{k}^{\prime}$ for $0\leq u\leq q+1;$%

\begin{equation}
t_{u}=\min\left\{  r\text{ }|\text{ }i_{r}+\ell_{\left(  r\right)  }%
-\ell_{\left(  o\left(  u\right)  \right)  }>{\epsilon}_{u}>i_{r}\right\}  ;
\label{ten}%
\end{equation}

$\hspace*{1.02in}o\left(  u\right)  =\max\left\{  r\text{ }|\text{ }i_{r}%
\geq\epsilon_{u}\right\}  $ and $o^{\prime}\left(  u\right)  =\max\left\{
r\text{ }|\text{ }\epsilon_{r}\leq i_{u}\right\}  .$\vspace{0.1in}\newline
Extend $\Delta$ multiplicatively to all of $C_{\ast}(K_{n+2})$, using the fact
that the cells of $K_{n+2}$ are products of cells $K_{i+2}$ with $i<n$.

Note that right-hand and left-hand factors in each component of $\Delta
T_{n+2}$ are expressed in first and second fundamental form, respectively. In
particular, the terms given by the extremes $p=0$ and $p=n$ are the primitive
terms of $\Delta T_{n+2}:$
\[
\left\{  d_{(0,1)}\cdots d_{(0,1)}\otimes1+1\otimes d_{(1,1)}\cdots
d_{(n,1)}\right\}  \left(  T_{n+2}\otimes T_{n+2}\right)  .
\]
The sign in $\Delta T_{n+2}=\sum(-1)^{\epsilon}b\otimes a$ is the product of
five signs: $\left(  -1\right)  ^{\epsilon}=sgn\left(  b\right)  \cdot
sgn\left(  b_{0}\right)  \cdot sgn(b_{0},a_{1})\cdot sgn\left(  a_{1}\right)
\cdot sgn\left(  a\right)  ,$ where the face $b_{0}$ is obtained from $b$ by
using the same $\varepsilon$'s but all $i_{k}^{\prime}=0$ (i.e., $\ell
_{k}^{\prime}$ is replaced by $\ell_{k}^{\prime}+i_{k}^{\prime}-i_{k-1}%
^{\prime}$)$,$ the face $a_{1}$ is obtained from $a$ by using the same $i_{r}%
$'s but all $\ell_{r}=1$ and $sgn(b_{0},a_{1})$ is the sign of the unshuffle
$\{i_{p}<\cdots<i_{1},$ $\epsilon_{1}<\epsilon_{2}<\cdots<\epsilon_{q}\}{.}$
Geometrically, $b_{0}$ and $a_{1}$ lie on orthogonal faces of the cube $I^{n}$
and are uniquely defined by the property that the canonical cellular
projection $K_{n+2}\rightarrow I^{n}$ maps $b_{0}\mapsto\left(  x_{1}%
,...,x_{\epsilon_{1}-1},0,...,x_{\epsilon_{q}-1},0,...,x_{n}\right)  $ and
$a_{1}\mapsto\left(  x_{1},...,x_{i_{p}-1},1,...,x_{i_{1}-1},1,...,x_{n}%
\right)  $.

\begin{example}
We obtain by direct calculation:

\noindent For $T_{3}\in C_{1}\left(  K_{3}\right) : $
\[
\Delta T_{3} =( d_{\left(  0,1\right) } \otimes1+1\otimes d_{\left(
1,1\right) })(T_{3}\otimes T_{3}).
\]

\noindent For $T_{4}\in C_{2}\left(  K_{4}\right) : $
\begin{align*}
\Delta T_{4}  &  =( d_{\left(  0,1\right)  }d_{\left(  0,1\right) }
\otimes1+1\otimes d_{\left(  1,1\right)  }d_{\left(  2,1\right) }+ d_{\left(
0,2\right) }\otimes d_{\left(  1,1\right) }\\
&  \qquad\qquad+d_{\left(  0,2\right)  }\otimes d_{\left(  1,2\right)
}^{\mathstrut}+d_{\left(  1,1\right)  }\otimes d_{\left(  1,2\right)
}^{\mathstrut}-d_{\left(  0,1\right)  }\otimes d_{\left(  2,1\right)
}^{\mathstrut} ) (T_{4}\otimes T_{4}).
\end{align*}
For $T_{5}\in C_{3}\left(  K_{5}\right) : $\vspace{-0.2in}

\[%
\begin{array}
[c]{llll}%
\Delta T_{5}= & ( d_{\left(  0,1\right) }d_{\left(  0,1\right)  }d_{\left(
0,1\right)  }\otimes1+1\otimes d_{\left(  1,1\right)  }d_{\left(  2,1\right)
}d_{\left(  3,1\right) }+ &  & \\
& +d_{\left(  0,3\right)  }\otimes d_{\left(  1,1\right)  }d_{\left(
2,2\right)  }-d_{\left(  1,1\right)  }\otimes d_{\left(  1,2\right)
}d_{\left(  3,1\right)  }+d_{\left(  0,1\right)  }d_{\left(  0,2\right)
}\otimes d_{\left(  1,2\right)  } &  & \\
& +d_{\left(  0,3\right)  }\otimes d_{\left(  1,2\right)  }d_{\left(
2,1\right)  }-d_{\left(  1,2\right)  }\otimes d_{\left(  1,2\right)
}d_{\left(  2,1\right)  }+d_{\left(  0,1\right)  }d_{\left(  0,2\right)
}\otimes d_{\left(  1,3\right)  } &  & \\
& +d_{\left(  0,3\right)  }\otimes d_{\left(  1,1\right)  }d_{\left(
2,1\right)  }-d_{\left(  1,2\right)  }\otimes d_{\left(  1,1\right)
}d_{\left(  2,2\right)  }-d_{\left(  0,2\right)  }d_{\left(  0,1\right)
}\otimes d_{\left(  2,1\right)  } &  & \\
& +d_{\left(  0,1\right)  }\otimes d_{\left(  2,1\right)  }d_{\left(
3,1\right)  }+d_{\left(  0,1\right)  }d_{\left(  0,1\right)  }\otimes
d_{\left(  3,1\right)  }+d_{\left(  0,2\right)  }d_{\left(  0,1\right)
}\otimes d_{\left(  2,2\right)  } &  & \\
& +d_{\left(  2,1\right)  }\otimes d_{\left(  1,1\right)  }d_{\left(
2,2\right)  }+d_{\left(  1,1\right)  }d_{\left(  0,1\right)  }\otimes
d_{\left(  2,2\right)  }+d_{\left(  0,2\right)  }d_{\left(  1,1\right)
}\otimes d_{\left(  1,2\right)  } &  & \\
& -d_{\left(  0,2\right)  }\otimes d_{\left(  1,1\right)  }d_{\left(
3,1\right)  }+d_{\left(  1,1\right)  }d_{\left(  1,1\right)  }\otimes
d_{\left(  1,3\right)  }+d_{\left(  0,2\right)  }d_{\left(  1,1\right)
}\otimes d_{\left(  1,3\right)  } &  & \\
& -d_{\left(  0,2\right)  }\otimes d_{\left(  1,2\right)  }d_{\left(
3,1\right)  }+d_{\left(  0,1\right)  }d_{\left(  0,2\right)  }\otimes
d_{\left(  1,1\right)  } )(T_{5}\otimes T_{5}). &  &
\end{array}
\]

\end{example}

$\vspace{0.1in}$

\noindent Our main result, stated as the following theorem, is proved in
Section \ref{pdelta}:

\begin{theorem}
\label{Delta} For each $n\geq0,$ the map $\Delta:C_{\ast}\left(
K_{n+2}\right)  \rightarrow C_{\ast}\left(  K_{n+2}\right)  \otimes C_{\ast
}\left(  K_{n+2}\right)  $ defined above is a chain map.
\end{theorem}

\noindent Identify the sequence of cellular chain complexes $\left\{  C_{\ast
}(K_{n})\right\}  _{n\geq2}$ with the $A_{\infty}$-operad $\mathcal{A}
_{\infty}$ \cite{May}. Since $\Delta$ is extended multiplicatively on
decomposable faces, we immediately obtain:

\begin{corollary}
\label{corollary1}The sequence of chain maps $\left\{  \Delta:C_{\ast}%
(K_{n})\rightarrow C_{\ast}(K_{n})\otimes C_{\ast}(K_{n})\right\}  _{n\geq2}$
induces a morphism of operads $\mathcal{A}_{\infty}\rightarrow\mathcal{A}%
_{\infty}\otimes\mathcal{A}_{\infty}.$
\end{corollary}

\section{A proof of Theorem \ref{Delta}}

\label{pdelta}

In this section we prove that the diagonal $\Delta$ defined in Section
\ref{diag} is a chain map. We begin with some preliminaries. It will be
convenient to rewrite face relation (3) as follows:%
\[
\qquad d_{\left(  i_{q+1},\ell_{q+1}\right)  }^{q+1}d_{\left(  i_{q},\ell
_{q}\right)  }^{q}=d_{\left(  i_{q}-\ell_{q+1},\ell_{q}\right)  }%
^{q+1}d_{\left(  i_{q+1},\ell_{q+1}\right)  }^{q};\qquad i_{q}>i_{q+1}%
+\ell_{q+1}.\qquad\left(  3^{\prime}\right)
\]

\begin{definition}
For $1\leq k\leq m,$ the $k^{th}$ l\underline{eft-transfer} $\tau_{\ell}^{k}$
of a composition $d_{(i_{m},\ell_{m})}\linebreak\cdots d_{\left(  i_{2}%
,\ell_{2}\right)  }d_{(i_{1},\ell_{1})}$ in first fundamental form is one of
the following compositions:\vspace{0.1in}

\begin{enumerate}
\item[(a)] If $i_{k+1}+\ell_{k+1}\geq i_{k}$\text{, apply face relation }%
$(2)$\text{\ to }$d^{k+1}d^{k}$\text{\ }and obtain
\[
\cdots d_{\left(  i_{k}-i_{k+1},\ell_{k}\right)  }^{k}d_{\left(  i_{k+1}%
,\ell_{k+1}+\ell_{k}\right)  }^{k}\cdots;
\]
then successively apply face relation\text{ }$(1)$\text{\ to }$d^{j+2}d^{j}%
$\text{\ for }$j=k,k+1,\ldots,m-2,$ and obtain
\[
\tau_{\ell}^{k}=d_{\left(  i_{k}-i_{k+1},\ell_{k}\right)  }^{k}d_{(i_{m}%
,\ell_{m})}^{m-1}\cdots d_{\left(  i_{k+2},\ell_{k+2}\right)  }^{k+1}%
d_{\left(  i_{k+1},\ell_{k+1}+\ell_{k}\right)  }^{k}\cdots d_{(i_{1},\ell
_{1})}^{1}.
\]

\item[(b)] If\text{ }$i_{k+1}+\ell_{k+1}<i_{k}\leq i_{q+1}+\ell_{\left(
q+1\right)  }-\ell_{\left(  k\right)  }$\text{\ for some smallest integer
}$q>k,$ successively apply face relation\text{ }$(3^{\prime})$\text{\ to
}$d^{j+1}d^{j}$\text{\ for }$j=k,k+1,\ldots,q,$\text{\ }and obtain
\[
\text{\hspace*{0.3in}}\cdots d_{(i_{q+2},\ell_{q+2})}^{q+2}d_{(i_{k}%
+\ell_{\left(  k\right)  }-\ell_{\left(  q+1\right)  },\ell_{k})}%
^{q+1}d_{(i_{q+1},\ell_{q+1})}^{q}\cdots d_{(i_{k+1},\ell_{k+1})}%
^{k}d_{(i_{k-1},\ell_{k-1})}\cdots.
\]
Apply face relation\text{ }$(2)$\text{\ to }$d^{q+1}d^{q};$\text{\ }then
successively apply face relation (1)\ to $d^{j+2}d^{j}$\text{ for
}$j=q,q+1,\ldots,m-2,$ and obtain
\[
\tau_{\ell}^{k}=d_{\left(  i,\ell_{k}\right)  }^{q}d_{(i_{m},\ell_{m})}%
^{m-1}\cdots d_{(i_{k+1},\ell_{k+1})}^{k}d_{(i_{k-1},\ell_{k-1})}\cdots
d_{(i_{1},\ell_{1})},
\]
where\text{ }$i=i_{k}-i_{q+1}+\ell_{\left(  k\right)  }-\ell_{\left(
q\right)  }.$\vspace{0.1in}

\item[(c)] Otherwise, successively apply face relation\text{ }$\left(
3^{\prime}\right)  $\text{\ to }$d^{j+1}d^{j}$\text{\ for }$j=k,k+1,\ldots
,m-1,$ and obtain
\[
\tau_{\ell}^{k}=d_{(i,\ell_{k})}^{m}d_{(i_{m},\ell_{m})}^{m-1}\cdots
d_{(i_{k+1},\ell_{k+1})}^{k}d_{(i_{k-1},\ell_{k-1})}\cdots d_{(i_{1},\ell
_{1})}%
\]
where\text{ }$i=i_{k}+\ell_{\left(  k\right)  }-\ell_{\left(  m\right)  }$.
\end{enumerate}
\end{definition}

\begin{definition}
For $1\leq k\leq m,$ the $k^{th}$ \underline{right-transfer} $\tau_{r}^{k}$ of
a composition $d_{(i_{m},\ell_{m})}\linebreak\cdots d_{\left(  i_{2},\ell
_{2}\right)  }d_{(i_{1},\ell_{1})}$ in first fundamental form is one of the
following compositions:\vspace{0.1in}

\begin{enumerate}
\item[(a)] If\text{ }$i_{p}+\ell_{\left(  p\right)  }\leq i_{k}+\ell_{\left(
k\right)  }<i_{p-1}+\ell_{\left(  p-1\right)  }$\text{\ for some greatest
integer }$1<p\leq k,$ successively apply face relation\text{ }$\left(
3^{\prime}\right)  $ to\text{ }$d^{j+1}d^{j}\ $for\text{ }$j=p-1,\ldots,2,1; $
$j=p,\ldots,3,2;$ $\ \ldots,$\text{\ }$\ j=k-1,k-2,\ldots,k-\left(
p-1\right)  $\text{. }Then apply face relation $(2)$\text{\ to }$d^{j+1}d^{j}%
$\text{\ for }$j=k-p,\ldots,2,1$ and obtain%
\[%
\begin{tabular}
[c]{l}%
$\tau_{r}^{k}=d_{(i_{m},\ell_{m})}\cdots d_{(i_{k+1},\ell_{k+1})}%
d_{(i_{p-1}-\ell,\ell_{p-1})}^{k}$\\
\\
$\qquad\qquad\qquad\qquad\cdots d_{(i_{1}-\ell,\ell_{1})}^{k-p+2}%
d_{(i_{k-1}-i_{k},\ell_{k-1})}^{k-p}\cdots d_{(i_{p}-i_{k},\ell_{p})}%
^{1}d_{(i_{k},\ell)}^{1},$%
\end{tabular}
\]
where\text{ }$\ell=\ell_{\left(  k\right)  }-\ell_{\left(  p-1\right)  }$.
\vspace{0.1in}

\item[(b)] \text{Otherwise, successively apply face relation }$(2)$\text{ \ to
}$d^{j+1}d^{j}$\text{\ for }$j=k-1,$\linebreak$k-2,$ $\ldots,2,1,$ \text{and
obtain}
\[
\tau_{r}^{k}=d_{(i_{m},\ell_{m})}\cdots d_{(i_{k+1},\ell_{k+1})}d_{\left(
i_{k-1}-i_{k},\ell_{k-1}\right)  }^{k-1}\cdots d_{(i_{1}-i_{k},\ell_{1})}%
^{1}d_{\left(  i_{k},\ell_{\left(  k\right)  }\right)  }^{1}.
\]
\vspace{0.05in}
\end{enumerate}
\end{definition}

Note that if $I=\left(  i_{1},\ell_{1}\right)  ,\ldots,\left(  i_{m},\ell
_{m}\right)  $ is a type I sequence and $a=d_{I}\left(  T_{n+2}\right)  \in
C_{n-m}\left(  K_{n+2}\right)  ,$ then for each $k,$ the $k^{th}$ right
transfer $\tau_{r}^{k}\left(  T_{n+2}\right)  $ expresses $a$ in first
fundamental form as a face of some $\left(  n-1\right)  $-face $b$ of
$T_{n+2}.$ The expressions $\tau_{r}^{1}\left(  T_{n+2}\right)  ,\ldots
,\tau_{r}^{m}\left(  T_{n+2}\right)  $ determine the $m$ distinct $\left(
n-1\right)  $-faces $b$ containing $a$.

There are analogous left and right transfers for compositions in second
fundamental form.

\begin{definition}
For $1\leq k\leq m,$ the $k^{th}$ \underline{left-transfer} $\tau_{\ell}^{k}$
of a composition $d_{(i_{m},\ell_{m})}\linebreak\cdots d_{\left(  i_{2}%
,\ell_{2}\right)  }d_{(i_{1},\ell_{1})}$ in second fundamental form is one of
the following compositions:\vspace{0.1in}

\begin{enumerate}
\item[(a)] If $i_{k+1}\leq i_{k},$\text{ apply face relation }$(2)$\text{\ to
}$d^{k+1}d^{k}$\text{, then successively apply face}\linebreak\text{relation
}$(1)$\text{\ to }$d^{j+2}d^{j}$\text{\ for }$j=k,k+1,\ldots,m-2,$ \text{and
obtain}
\[
\tau_{\ell}^{k}=d_{\left(  i_{k}-i_{k+1},\ell_{k}\right)  }^{k}d_{(i_{m}%
,\ell_{m})}^{m-1}\cdots d_{\left(  i_{k+2},\ell_{k+2}\right)  }^{k+1}%
d_{\left(  i_{k+1},\ell_{k+1}+\ell_{k}\right)  }^{k}\cdots d_{(i_{1},\ell
_{1})}^{1}.
\]

\item[(b)] If $i_{k+1}>i_{k}\geq i_{q+1}$ \text{for some smallest integer
}$q>k,$ \text{successively apply face}\linebreak\text{relation }%
$(3)$\text{\ to }$d^{j+1}d^{j}$\text{\ for }$j=k,k+1,\ldots,q.$ \text{Apply
face relation }$(2)$\text{\ to}\linebreak$d^{q+1}d^{q};$\text{\ then
successively apply face relation }$(1)$\text{\ to }$d^{j+2}d^{j}$ \text{\ for
}$j=q,q+1,\ldots,m-2,$ \text{and obtain}%
\[%
\begin{tabular}
[c]{l}%
$\hspace*{0.3in}\tau_{\ell}^{k}=d_{\left(  i_{k}-i_{q+1},\ell_{k}\right)
}^{q}d_{(i_{m},\ell_{m})}^{m-1}\cdots d_{\left(  i_{q+2},\ell_{q+2}\right)
}^{q+1}d_{(i_{q+1},\ell_{q+1}+\ell_{k})}^{q}d_{(i_{q}+\ell_{k},\ell_{q}%
)}^{q-1}$\\
\\
$\hspace*{2in}\cdots d_{(i_{k+1}+\ell_{k},\ell_{k+1})}^{k}d_{(i_{k-1}%
,\ell_{k-1})}\cdots d_{(i_{1},\ell_{1})}.$%
\end{tabular}
\]
\vspace{0.1in}

\item[(c)] \text{Otherwise}, \text{successively apply face relation }$\left(
3\right)  $\text{\ to }$d^{j+1}d^{j}$\text{\ for }$j=k,k+1,$\linebreak%
$\ldots,m-1,$ \text{and obtain}%
\[%
\begin{tabular}
[c]{l}%
$\hspace*{0.38in}\tau_{\ell}^{k}=d_{(i_{k},\ell_{k})}^{m}d_{(i_{m}+\ell
_{k},\ell_{m})}^{m-1}\cdots d_{(i_{k+2}+\ell_{k},\ell_{k+2})}^{k+1}%
d_{(i_{k+1}+\ell_{k},\ell_{k+1})}^{k}d_{(i_{k-1},\ell_{k-1})}\cdots
d_{(i_{1},\ell_{1})}.$%
\end{tabular}
\
\]

\end{enumerate}
\end{definition}

\begin{definition}
For $1\leq k\leq m,$ the $k^{th}$ \underline{right-transfer} $\tau_{r}^{k}$ of
a composition $d_{(i_{m},\ell_{m})}\linebreak\cdots d_{\left(  i_{2},\ell
_{2}\right)  }d_{(i_{1},\ell_{1})}$ in second fundamental form is one of the
following compositions:\vspace{0.1in}

\begin{enumerate}
\item[(a)] \text{If }$i_{p-1}<i_{k}\leq i_{p}$\text{\ for some greatest
integer }$1<p\leq k,$ \text{successively apply}\linebreak\text{face relation
}$\left(  3\right)  $\text{\ to }$d^{j+1}d^{j}\ $\text{ for }$j=p-1,\ldots
,2,1;$ $j=p,\ldots,3,2;$\text{\ }$\ \ldots;$\text{ }$\ j=k-1,k-2,\ldots
,k-\left(  p-1\right)  $\text{. Then apply face relation }$(2) $ \text{\ to
}$d^{j+1}d^{j}$ \text{for }$j=k-p,\ldots,2,1$ \text{and obtain}%
\[%
\begin{tabular}
[c]{l}%
$\hspace*{0.3in}\tau_{r}^{k}=d_{(i_{m},\ell_{m})}\cdots d_{(i_{k+1},\ell
_{k+1})}d_{(i_{p-1},\ell_{p-1})}^{k}\cdots d_{(i_{1},\ell_{1})}^{k-p+2}%
d_{(i_{k-1}-i_{k},\ell_{k-1})}^{k-p}$\\
\\
$\hspace*{2in}\cdots d_{(i_{p}-i_{k},\ell_{p})}^{1}d_{(i_{k}+\ell_{\left(
p-1\right)  },\ell_{\left(  k\right)  }-\ell_{\left(  p-1\right)  })}^{1}.$%
\end{tabular}
\]
\vspace{0.1in}

\item[(b)] \text{Otherwise, successively apply face relation }$(2)$\text{ \ to
}$d^{j+1}d^{j}$\text{\ for }$j=k-1,\ldots,$\linebreak$2,1,$ \text{and
obtain}:
\[
\tau_{r}^{k}=d_{(i_{m},\ell_{m})}\cdots d_{(i_{k+1},\ell_{k+1})}%
d_{(i_{k-1}-i_{k},\ell_{k-1})}^{k}\cdots d_{(i_{1}-i_{k},\ell_{1})}%
^{1}d_{(i_{k},\ell_{\left(  k\right)  })}^{1}.
\]

\end{enumerate}
\end{definition}

Once again, if $I=\left(  i_{1},\ell_{1}\right)  ,\ldots,\left(  i_{m}%
,\ell_{m}\right)  $ is a type II sequence and $a=d_{I}\left(  T_{n+2}\right)
\in C_{n-m}\left(  K_{n+2}\right)  ,$ then for each $k,$ the $k^{th}$
right-transfer $\tau_{r}^{k}\left(  T_{n+2}\right)  $ expresses $a$ in second
fundamental form as a face of some $\left(  n-1\right)  $-face $b;$ the
expressions $\tau_{r}^{1}\left(  T_{n+2}\right)  ,\ldots,\linebreak\tau
_{r}^{m}\left(  T_{n+2}\right)  $ determine the $m$ distinct $\left(
n-1\right)  $-faces $b$ containing $a$. Thus two $\left(  n-m\right)  $-faces
$a_{1}$ and $a_{2}$ expressed in first and second fundamental form,
respectively, are contained in the same $\left(  n-1\right)  $-face $b$ if and
only if there exist right transfers $\tau_{r}^{k_{j}}$ such that $a_{j}%
=\tau_{r}^{k_{j}}\left(  T_{n+2}\right)  ,$ $j=1,2,$ and $b=d_{(i,\ell
)}(T_{n+2}),$ where $d_{(i,\ell)}$ is the right-most face operator common to
$\tau_{r}^{k_{1}}$ and $\tau_{r}^{k_{2}}.$ We state this formally in the
following lemma:

\begin{lemma}
\label{l1}Let $0\leq m_{1},m_{2}\leq n$ and assume that $a_{1}=d_{(i_{m_{1}%
},\ell_{m_{1}})}\cdots d_{(i_{1},\ell_{1})}\left(  T_{n+2}\right)  $ and
$a_{2}=d_{(i_{m_{2}}^{\prime},\ell_{m_{2}}^{\prime})}\cdots d_{(i_{1}^{\prime
},\ell_{1}^{\prime})}\left(  T_{n+2}\right)  $ are expressed in first and
second fundamental form, respectively. Then $a_{1}$ and $a_{2}$ are contained
in the same $(n-1)$-face\ $d_{(i,\ell)}(T_{n+2})$ if and only if there exist
integers $k_{j}\leq m_{j}$ and greatest integers $1\leq p_{j}<k_{j}$ for
$j=1,2$, such that
\[%
\begin{array}
[c]{c}%
i_{p_{2}}^{\prime}<i_{k_{2}}^{\prime},\\
i_{k_{1}}+\ell_{\left(  k_{1}\right)  }-\ell_{\left(  p_{1}\right)  }%
<i_{p_{1}}%
\end{array}
\]
and
\[%
\begin{array}
[c]{l}%
i=i_{k_{1}}=i_{k_{2}}^{\prime}+\ell_{\left(  p_{2}\right)  }^{\prime},\\
\ell=\ell_{\left(  k_{1}\right)  }-\ell_{\left(  p_{1}\right)  }=\ell_{\left(
k_{2}\right)  }^{\prime}-\ell_{\left(  p_{2}\right)  }^{\prime}.
\end{array}
\]

\end{lemma}

Consider a solution $((i_{1},\ell_{1}),\ldots,(i_{p},\ell_{p});(i_{1}^{\prime
},\ell_{1}^{\prime}),\ldots,(i_{q}^{\prime},\ell_{q}^{\prime}))$ of system
(\ref{shuffle2}) and its related $\left(  p,q\right)  $-unshuffle
$\{i_{p}<\cdots<i_{1},$ $\epsilon_{1}<\epsilon_{2}<\cdots<\epsilon_{q}\}.$
Given $1\leq k\leq q+1,$ let $k_{1}=k-1,$ $k_{j+1}=o^{\prime}\left(  t_{k_{j}%
}\right)  $ for $j\geq1$ and note that $k_{j+1}<k_{j}$ for all $j.$ For $0\leq
m<\ell_{k}^{\prime},$ consider the following \textit{selection algorithm}:%
\[%
\begin{tabular}
[c]{l}%
\textbf{if}\textit{ }$k=1$\textit{ }\textbf{then}\textit{ }$z=p+1-i_{1}%
^{\prime}-m$\\
\textit{\hspace*{0.15in}}\textbf{else }$z=n+2;$ $j=0$\\
\textit{\hspace*{0.15in}\hspace*{0.15in}}\textbf{repeat}\\
\textit{\hspace*{0.15in}\hspace*{0.15in}\hspace*{0.15in}}$j\leftarrow
j+1\vspace{0.05in}$\\
\textit{\hspace*{0.15in}\hspace*{0.15in}\hspace*{0.15in}}\textbf{if}\textit{
}$i_{k_{j}}^{\prime}<i_{k}^{\prime}+m$\textit{ }\textbf{then}\textit{\ }%
$i_{z}=\epsilon_{k_{j}}-i_{k_{j}}^{\prime}+i_{k}^{\prime}+m\vspace{0.05in}$\\
\textit{\hspace*{0.15in}\hspace*{0.15in}\hspace*{0.15in}}\textbf{if}\textit{
}$i_{t_{k_{j}}}-\ell_{\left(  k_{j+1}\right)  }^{\prime}=i_{k}^{\prime}%
+m$\textit{ }\textbf{then }$z=t_{k_{j}}\vspace{0.05in}$\\
\textit{\hspace*{0.15in}\hspace*{0.15in}}\textbf{until}\textit{ }$z<n+2$\\
\textbf{endif}%
\end{tabular}
\ \
\]
It will be clear from the proof of Lemma \ref{l2} below that the selection
algorithm eventually terminates.

\begin{example}
Let $p=q=4$ and consider the following solution of system (\ref{shuffle2}):%
\[
\left(  \left(  7,1\right)  ,\left(  6,1\right)  ,\left(  4,2\right)  ,\left(
2,3\right)  ;\left(  0,1\right)  ,\left(  1,1\right)  ,\left(  1,2\right)
,\left(  0,4\right)  \right)  .
\]
Then $t_{1}=5,$ $t_{2}=4,$ $t_{3}=3,$ $t_{4}=4$ and the \textit{selection
algorithm} produces the following:\vspace{0.1in}%
\[%
\begin{tabular}
[c]{|l||l|l|l|l|l|l|l|l|}\hline
$k$ & 1 & 2 & 3 & 3 & 4 & 4 & 4 & 4\\\hline
$m$ & 0 & 0 & 0 & 1 & 0 & 1 & 2 & 3\\\hline
$z$ & 5 & 5 & 4 & 3 & 5 & 4 & 2 & 1\\\hline
\end{tabular}
\
\]
\vspace{0.1in}
\end{example}

The key to our proof that $\Delta$ is a chain map is given by our next lemma.

\begin{lemma}
\label{l2}Given $1\leq k\leq q+1$ and $0\leq m<\ell_{k}^{\prime},$ let $z$ be
the integer given by the selection algorithm.

\begin{enumerate}
\item[(a)] $z>o(k)$\vspace{0.1in}

\item[(b)] $i_{z}+\ell_{(z)}-\ell_{(o(k))}\geq\epsilon_{k}.$\vspace{0.1in}

\item[(c)] $i_{k}^{\prime}+m=i_{z}-\ell_{(o^{\prime}(z))}^{\prime};$%
\vspace{0.1in}

\item[(d)] If $o^{\prime}(z)<r<k,$ then $i_{k}^{\prime}+m\leq i_{r}^{\prime}%
;$\vspace{0.1in}

\item[(e)] If $i_{k}^{\prime}+m>\min\limits_{o^{\prime}\left(  t_{k}\right)
<r<k}\left\{  i_{r}^{\prime},\text{ }i_{t_{k}}-\ell_{(o^{\prime}\left(
t_{k}\right)  )}^{\prime}\right\}  ,$ then\vspace{0.1in}

\begin{enumerate}
\item[(e.1)] $z<t_{k}$ and $i_{z}+\ell_{(z)}-\ell_{(o\left(  k\right)
)}=\epsilon_{k};$\vspace{0.1in}

\item[(e.2)] $o\left(  k\right)  =\max\limits_{r<z}\left\{  r\text{
}|\text{\ }i_{r}>i_{z}+\ell_{(z)}-\ell_{(r)}\right\}  ;$\vspace{0.1in}

\item[(e.3)] $o^{\prime}(z)=\max\limits_{r<k}\left\{  r\text{ }|\text{ }%
i_{r}^{\prime}<i_{k}^{\prime}+m\right\}  .$\vspace{0.1in}
\end{enumerate}
\end{enumerate}
\end{lemma}

\noindent\underline{\textit{Proof:}}\vspace{0.1in}

(a)\textit{ }When $k=1,$ $z=\left(  p+1-\epsilon_{1}\right)  +\left(  \ell
_{1}^{\prime}-m\right)  >p+1-\epsilon_{1}=o\left(  1\right)  .$ For $k>1,$
first note that $i_{t_{k_{j}}}<\epsilon_{k_{j}}<\epsilon_{k}\leq i_{o\left(
k\right)  }$ for all $j\geq1,$ by the definition of $t_{k_{j}}$ in
(\ref{ten})$.$ So the result follows whenever $i_{z}\leq i_{t_{k_{j}}}$ for
some $j\geq1.$ If $i_{k-1}^{\prime}<i_{k}^{\prime}+m,$ then $i_{z}%
=\epsilon_{k-1}+i_{k}^{\prime}+m-i_{k-1}^{\prime}<\epsilon_{k-1}+i_{k}%
^{\prime}+\ell_{k}^{\prime}-i_{k-1}^{\prime}=\epsilon_{k}$. If $i_{k_{j}%
}^{\prime}<i_{k}^{\prime}+m$ for some $j\geq2$ and $i_{k}^{\prime}+m\leq
i_{k_{r}}^{\prime}\leq i_{t_{k_{r}}}-\ell_{\left(  k_{r+1}\right)  }^{\prime}$
for all $r<j,$ where the later follows from inequality (3) of system
(\ref{shuffle2}), then $i_{z}=\epsilon_{k_{j}}+i_{k}^{\prime}+m-i_{k_{j}%
}^{\prime}<\epsilon_{k_{j}}+i_{t_{k_{j-1}}}-\ell_{\left(  k_{j}\right)
}^{\prime}-i_{k_{j}}^{\prime}=i_{t_{k_{j-1}}}.$

(b) When $k=1,$ indices $0=i_{p+1}<i_{p}<\cdots<i_{p+2-\epsilon_{1}}$ are
consecutive; hence $i_{r}+\ell_{\left(  r\right)  }=i_{r-1}+\ell_{\left(
r-1\right)  }+\ell_{r}-1\geq i_{r-1}+\ell_{\left(  r-1\right)  }$ for
$p+3-\epsilon_{1}\leq r\leq p+1.$ Since $z\geq p+2-\epsilon_{1}$ we have
$i_{z}+\ell_{(z)}-\ell_{(o(1))}\geq i_{p+2-\epsilon_{1}}+\ell_{\left(
p+2-\epsilon_{1}\right)  }-\ell_{\left(  p+1-\epsilon_{1}\right)  }%
=\epsilon_{1}-1+\ell_{p+2-\epsilon_{1}}\geq\epsilon_{1}.$ So consider $k>1.$
If $i_{k-1}^{\prime}<i_{k}^{\prime}+m,$ then $\epsilon_{k-1}<i_{z}%
<\epsilon_{k}<i_{o\left(  k\right)  }$ so that $\epsilon_{k}=i_{z}+z-o\left(
k\right)  \leq i_{z}+\ell_{(z)}-\ell_{(o(k))}.$ If $i_{k-1}^{\prime}\geq
i_{k}^{\prime}+m$ and $z=t_{k-1},$ then $o\left(  k-1\right)  -o\left(
k\right)  =\epsilon_{k}-\epsilon_{k-1}-1$ so that $i_{z}+\ell_{(z)}%
-\ell_{(o(k))}=i_{z}+\ell_{(z)}-\ell_{(o(k-1))}+\ell_{(o(k-1))}-\ell
_{(o\left(  k\right)  )}>\epsilon_{k-1}+\epsilon_{k}-\epsilon_{k-1}%
-1=\epsilon_{k}-1$ by the definition of $t_{k-1}$ in (\ref{ten}). If
$i_{k_{2}}^{\prime}<i_{k}^{\prime}+m$ then $\epsilon_{k_{2}}<i_{z}%
<i_{t_{k_{1}}}<\epsilon_{k_{1}}$ so that $\ell_{\left(  z\right)  }%
-\ell_{\left(  t_{k_{1}}\right)  }\geq z-t_{k_{1}}=i_{t_{k_{1}}}-i_{z},$ where
equality holds iff $\ell_{t_{k_{1}}+1}=\cdots=\ell_{z}=1;$ hence $i_{z}%
+\ell_{\left(  z\right)  }-\ell_{\left(  o\left(  k\right)  \right)  }%
=i_{z}+\ell_{\left(  z\right)  }-\ell_{\left(  t_{k_{1}}\right)  }%
+\ell_{\left(  t_{k_{1}}\right)  }-\ell_{\left(  o\left(  k\right)  \right)
}\geq i_{z}+i_{t_{k_{1}}}-i_{z}+\ell_{\left(  t_{k_{1}}\right)  }%
-\ell_{\left(  o\left(  k\right)  \right)  }>\epsilon_{k}-1\ $by the previous
calculation. If $i_{k_{2}}^{\prime}\geq i_{k}^{\prime}+m$ and $z=t_{k_{2}},$
then $i_{z}+\ell_{\left(  z\right)  }-\ell_{\left(  o\left(  k_{1}\right)
\right)  }=i_{z}+\ell_{\left(  z\right)  }-\ell_{\left(  o\left(
k_{2}\right)  \right)  }+\ell_{\left(  o\left(  k_{2}\right)  \right)  }%
-\ell_{\left(  o\left(  k_{1}\right)  \right)  }>\epsilon_{k_{2}}%
+\ell_{\left(  o\left(  k_{2}\right)  \right)  }-\ell_{\left(  o\left(
k_{1}\right)  \right)  }\geq i_{t_{k_{1}}}+\ell_{\left(  t_{k_{1}}\right)
}-\ell_{\left(  o\left(  k_{2}\right)  \right)  }+\ell_{\left(  o\left(
k_{2}\right)  \right)  }-\ell_{\left(  o\left(  k_{1}\right)  \right)
}>\epsilon_{k_{1}}$ since $t_{k_{1}}<t_{k_{2}}$ implies $i_{t_{k_{1}}}%
+\ell_{\left(  t_{k_{1}}\right)  }-\ell_{\left(  o\left(  k_{2}\right)
\right)  }\leq\epsilon_{k_{2}}$ by the definition of $t_{k_{2}}.$ Thus
$i_{z}+\ell_{\left(  z\right)  }-\ell_{\left(  o\left(  k\right)  \right)
}=i_{z}+\ell_{\left(  z\right)  }-\ell_{\left(  o\left(  k_{1}\right)
\right)  }+\ell_{\left(  o\left(  k_{1}\right)  \right)  }-\ell_{\left(
o\left(  k\right)  \right)  }>\epsilon_{k}-1.$ Note that in general,
$i_{t_{k_{j}}}+\ell_{\left(  t_{k_{j}}\right)  }-\ell_{\left(  o\left(
k_{1}\right)  \right)  }>\epsilon_{k_{1}}$ so that $i_{t_{k_{j}}}%
+\ell_{\left(  t_{k_{j}}\right)  }-\ell_{\left(  o\left(  k\right)  \right)
}\geq\epsilon_{k}$ for all $j\geq1.$ So if necessary, continue in like manner
until the desired $z$ is found at which point the result follows.

(c) This result follows immediately from the choice of $z.$

(d) Note that $k>1$. If $i_{k_{1}}^{\prime}<i_{k}^{\prime}+m$, then
$o^{\prime}\left(  z\right)  =k-1$ and the case is vacuous. So assume that
$i_{k}^{\prime}+m\leq i_{k_{1}}^{\prime}.$ If either $i_{k}^{\prime
}+m=i_{t_{k_{1}}}-\ell_{\left(  k_{2}\right)  }^{\prime}$ or $i_{k_{2}%
}^{\prime}<i_{k}^{\prime}+m,$ then $o^{\prime}\left(  z\right)  =o^{\prime
}\left(  t_{k_{1}}\right)  $ and $i_{k_{1}}^{\prime}\leq i_{r}^{\prime}$ for
$o^{\prime}\left(  z\right)  <r<k.$ Otherwise, $i_{k}^{\prime}+m\leq
\min\left\{  i_{k_{1}}^{\prime},\text{ }i_{k_{2}}^{\prime}\right\}  .$ If
either $i_{k}^{\prime}+m=i_{t_{k_{2}}}-\ell_{\left(  k_{3}\right)  }^{\prime}$
or $i_{k_{3}}^{\prime}<i_{k}^{\prime}+m,$ then $o^{\prime}\left(  z\right)
=o^{\prime}\left(  t_{k_{2}}\right)  $ and $i_{k_{2}}^{\prime}\leq
i_{r}^{\prime}$ for $o^{\prime}\left(  z\right)  <r\leq k_{2}.$ But $i_{k_{1}%
}^{\prime}\leq i_{r}^{\prime}$ for $k_{2}=o^{\prime}\left(  t_{k_{1}}\right)
<r<k,$ so the desired inequality holds for $o^{\prime}\left(  z\right)  <r<k.$
Continue in this manner until $i_{k}^{\prime}+m\leq\min\left\{  i_{k_{1}%
}^{\prime},\ldots,i_{k_{j}}^{\prime}\right\}  $ for some $j,$ at which point
the conclusion follows.

(e) If $k=1,$ $o^{\prime}\left(  z\right)  =0$ so that $i_{1}^{\prime}%
+m=i_{z}$ and $i_{z}+\ell_{\left(  z\right)  }\geq\epsilon_{1}.$ Now
$i_{1}^{\prime}+m>i_{t_{1}}\ $by assumption, hence $z<t_{1}.$ But $t_{1}$ is
the smallest integer $r$ such that $i_{r}+\ell_{\left(  r\right)  }%
>\epsilon_{1};$ therefore $i_{z}+\ell_{\left(  z\right)  }=\epsilon_{1},$ by
(b). Let $k>1$ and suppose that $z\geq t_{k}.$ Then $i_{z}\leq i_{t_{k}},$ in
which case $o^{\prime}\left(  z\right)  =o^{\prime}\left(  t_{k}\right)  .$
But by (c), $i_{k}^{\prime}+m=i_{z}-\ell_{\left(  o^{\prime}\left(  z\right)
\right)  }^{\prime}\leq i_{t_{k}}-\ell_{\left(  o^{\prime}\left(
t_{k}\right)  \right)  }^{\prime}$ and by (d), $i_{k}^{\prime}+m\leq
i_{r}^{\prime}$ for $o^{\prime}\left(  t_{k}\right)  <r<k,$ contradicting the
hypothesis. Hence $z<t_{k}.$ But $t_{k}$ is the smallest possible integer such
that $i_{t_{k}}+\ell_{\left(  t_{k}\right)  }-\ell_{\left(  o\left(  k\right)
\right)  }>\epsilon_{k},$ so (e.1) follows from by (b). Results (e.2) and
(e.3) are obvious.\bigskip

\noindent\underline{\textit{Proof of Theorem \ref{Delta}:}}$\vspace{0.1in}$

Let $n\geq0.$ We show that if $u\otimes v$ is a component of $d^{\otimes
2}\Delta(T_{n+2})$ or $\Delta d(T_{n+2}),$ there is a corresponding component
that cancels it, in which case $(d^{\otimes2}\Delta-\Delta d)(T_{n+2})$ $=0.$
Let $a^{\prime}\otimes a$ be a component of $\Delta T_{n+2}.$ We consider
various cases.

\underline{\textit{Case I.}} Consider a component $d_{(i,\ell)}^{r}\otimes1$
of $d\otimes1,$ where $i+\ell\leq\ell_{r}^{\prime},\ i\geq0,$ $1\leq\ell
<\ell_{r}^{\prime}$\ and $1\leq r\leq q+1.$ Reduce $b^{\prime}=d_{(i,\ell
)}^{r}(a^{\prime})$ to an expression in second fundamental form, i.e., if
$1\leq r\leq q,$ successively apply relation (1) to $d^{j+1}d^{j}$ for
$\ j=q,...,r+1;$ apply (2) to replace $d_{(i,\ell)}^{r}d_{(i_{r}^{\prime}%
,\ell_{r}^{\prime})}^{r}$ by $d_{(i_{r}^{\prime},\ell_{r}^{\prime}-\ell
)}^{r+1}d_{(i_{r}^{\prime}+i,\ell)}^{r}$; then successively apply (3) to
$d^{j+1}d^{j}\ $for either $j=r-1,...,\beta$ if $i+\ell<\ell_{r}^{\prime}$ or
for $\ j=q,...,\beta$ if $r=q+1,$ where $\beta$ is the greatest integer such
that $1\leq\beta\leq r$ and $i+\ell\geq i_{\beta-1}^{\prime}$. Then for
$i+\ell<\ell_{r}^{\prime},$\ $1\leq r\leq q+1$ we have
\begin{align*}
b^{\prime}  &  =d_{(i_{q}^{\prime},\ell_{q}^{\prime})}^{q+1}\cdots
d_{(i_{r+1}^{\prime},\ell_{r+1}^{\prime})}^{r+2}d_{(i_{r}^{\prime},\ell
_{r}^{\prime}-\ell)}^{r+1}d_{(i_{r-1}^{\prime}-\ell,\ell_{r-1}^{\prime
})_{\mathstrut}}^{r}\\
&  \hspace*{1in}\cdots d_{(i_{\beta}^{\prime}-\ell,\ell_{\beta}^{\prime}%
)}^{\beta+1}d_{(i_{r}^{\prime}+i,\ell)}^{\beta}d_{(i_{\beta-1}^{\prime}%
,\ell_{\beta-1}^{\prime})}^{{\beta}-1}\cdots d_{(i_{1}^{\prime},\ell
_{1}^{\prime})}^{1}(T_{n+2}),
\end{align*}
and for $i+\ell=\ell_{r}^{\prime},$\ $1\leq r\leq q+1$ we have
\[
b^{\prime}=d_{(i_{q}^{\prime},\ell_{q}^{\prime})}^{q+1}\cdots d_{(i_{r+1}%
^{\prime},\ell_{r+1}^{\prime})}^{r+2}d_{(i_{r}^{\prime},\ell_{r}^{\prime}%
-\ell)}^{r+1}d_{(i_{r}^{\prime}+i,\ell)}^{r}d_{(i_{{r}-1}^{\prime},\ell
_{r-1}^{\prime})}^{{r}-1}\cdots d_{(i_{1}^{\prime},\ell_{1}^{\prime})}%
^{1}(T_{n+2}).
\]

\underline{\textit{Case Ia.}} Let $i+\ell<\ell_{r}^{\prime}$ and
$i_{r}^{\prime}+i+\ell>i_{\beta-1}^{\prime}.$ As in the case of $a^{\prime}$
above, the inequalities for lower indices in an expression of $b^{\prime}$ as
a face of $T_{n+2}$ in first fundamental form are strict, but whose number now
is increased by one. Obviously we will have that
\[
\epsilon_{\beta}=i_{r}^{\prime}+i+\ell_{(\beta-1)}^{\prime}+\ell=i_{k}%
\]
for certain $1\leq k\leq p.$ Apply the $k^{th}$ left-transfer to obtain
\[
a=\tau_{\ell}^{k}(T_{n+2})=d_{(\tilde{\imath},\tilde{\ell})}^{\tilde{k}}(b),
\]
where $\tilde{k}=\alpha-1,$ $\tilde{\imath}=i_{k}-i_{\alpha}-\ell_{(\alpha
-1)}-\ell_{(k)},$ $\tilde{\ell}=\ell_{k}$ and\ $\alpha$ is the smallest
integer $k\leq\alpha\leq p+1$ with $i_{\alpha}+\ell_{(\alpha)}-\ell_{(k)}\geq
i_{k}.$ Then $b^{\prime}\otimes b$ is a component of $\Delta T_{n+2}$ as well.

\underline{\textit{Case Ib.}} Let $i+\ell<\ell_{r}^{\prime}$ and
$i_{r}^{\prime}+i+\ell=i_{\beta-1}^{\prime}.$ Apply the $(\beta-1)^{th}$
left-transfer to obtain
\[
b^{\prime}=\tau_{\ell}^{\beta-1}(T_{n+2})=d_{(\tilde{\imath},\tilde{\ell}%
)}^{\beta-1}(c^{\prime}),
\]
where $\tilde{\imath}=\ell$ and $\tilde{\ell}={\ell^{\prime}}_{\beta-1}.$ Then
$c^{\prime}\otimes a$ is also a component of $\Delta T_{n+2}$.

\underline{\textit{Case Ic.}} Let $i+\ell=\ell_{r}^{\prime};$ there are two subcases:

\underline{\textit{Subcase i.}} If $1\leq r\leq q$ and $i_{r}^{\prime}%
+i\leq\min(i_{t_{r}}-\ell_{(o^{\prime}(t_{r}))}^{\prime},i_{j}^{\prime
}),\ o^{\prime}(t_{r})<j<r,$ apply the $(r+1)^{th}$ left-transfer to obtain
\[
b^{\prime}=\tau_{\ell}^{r+1}(T_{n+2})=d_{(\tilde{\imath},\tilde{\ell}%
)}^{\tilde{r}}(c^{\prime}),
\]
where $\tilde{r}=\alpha,\ \tilde{\imath}=i_{r}^{\prime}-i_{\alpha}^{\prime
},\ \tilde{\ell}=\ell_{r}^{\prime}-\ell$ and\ $\alpha$ is the smallest integer
$r\leq\alpha\leq p+1$ with $i_{\alpha}^{\prime}\leq i_{r}^{\prime}.$ Although
certain $i^{\prime}$'s are increased by $i,$ while $\ell_{r}^{\prime}$ is
reduced by $i$ to obtain $\ell,$ it is straightforward to check that
$c^{\prime}\otimes a$ is also a component of $\Delta T_{n+2}$.

\underline{\textit{Subcase ii.}} Suppose that $1\leq r\leq q+1$ and
$i_{r}^{\prime}+i>\min(i_{t_{r}}-\ell_{(o^{\prime}(t_{r}))}^{\prime}%
,i_{j}^{\prime}),\ o^{\prime}(t_{r})<j<r.$ In view of Lemma \ref{l2} (with
$k=r$ and $m=i$) we have
\[
\epsilon_{r}=i_{z}+\ell_{(z)}-\ell_{(o(r))}\text{ and }i_{r}^{\prime}%
+i=i_{z}-\ell_{(o^{\prime}(z))}^{\prime},
\]
from which we also establish the equality%
\[
\ell_{(z)}-\ell_{(o(r))}=\ell_{(r-1)}^{\prime}+\ell-\ell_{(o^{\prime}%
(z))}^{\prime}.
\]
The hypotheses of Lemma \ref{l1} are satisfied by setting $k_{1}=z,\ p_{{1}%
}=o(r)$ and $k_{2}=r,\ p_{{2}}=o^{\prime}(z)$. Hence, $b^{\prime}\otimes a$ is
a component of $\Delta d_{(\tilde{\imath},\tilde{\ell})}(T_{n+2})$ with
$\tilde{\imath}=i_{z}$ and $\tilde{\ell}=\ell_{(z)}-\ell_{(o(r))}.$

\underline{\textit{Case II.}} Consider a component $1\otimes d_{(i,\ell)}^{r}$
of $1\otimes d,$ where $i+\ell\leq\ell_{r},\ i\geq0,$ $1\leq\ell<\ell_{r}%
,$\ and $1\leq r\leq p+1.$ Reduce $b=d_{(i,\ell)}^{r}(a)$ to the first
fundamental form, i.e., if $1\leq r\leq p,$ successively apply relation (1) to
$d^{j+1}d^{j}$ for $j=p,...,r+1,$ and apply (2) to replace $d_{(i,\ell)}%
^{r}d_{(i_{r},\ell_{r})}^{r}$ by $d_{(i_{r},\ell_{r}-\ell)}^{r+1}%
d_{(i_{r}+i,\ell)}^{r}.$ Then if $i>0,$ successively apply (3) to
$d^{j+1}d^{j},\ $for $j=r-1,...,\beta;$ or if $r=p+1,$ successively apply (3)
to $d^{j+1}d^{j}\ $for $j=p,...,\beta,$ where $\beta$ is the greatest integer
with $1\leq\beta\leq r$ and $i_{r}+i+\ell_{r-1}+\cdots+\ell_{\beta}\leq
i_{\beta-1},$ that is $\beta$ is determined by (\ref{ksign2}). Then for
$\beta=r,$ we have $i_{r}+i\leq i_{r-1}$; and for $i>0,$\ $1\leq r\leq p+1,$
\begin{align*}
b  &  =d_{(i_{p},\ell_{p})}^{p+1}\cdots d_{(i_{r+1},\ell_{r+1})}%
^{r+2}d_{(i_{r},\ell_{r}-\ell)}^{r+1}d_{(i_{r-1},\ell_{r-1})_{\mathstrut}}%
^{r}\\
&  \hspace*{1in}\cdots d_{(i_{{\beta}},\ell_{\beta})}^{\beta+1}d_{(\tilde{
\imath},{\ell})}^{\beta}d_{(i_{{\beta}-1},\ell_{\beta-1})}^{\beta-1}\cdots
d_{(i_{1},\ell_{1})}^{1}(T_{n+2}),
\end{align*}
where $\tilde{\imath}=i_{r}+i+\ell_{(r-1)}-\ell_{(\beta-1)},$ and
\[
b=d_{(i_{p-1},\ell_{p-1})}^{p}\cdots d_{(i_{{r}+1},\ell_{{r}+1})}^{{r}%
+2}d_{(i_{{r}},\ell_{{r}}-\ell)}^{{r+1}}d_{({i}_{r},{\ell})}^{{r}}d_{(i_{{
r}-1},\ell_{{r}-1})}^{{r}-1}\cdots d_{(i_{1},\ell_{1})}^{1}(T_{n+2}),
\]
for $i=0,\ 1\leq r\leq p+1.$

\underline{\textit{Case IIa.}} Let $i>0$ and $i_{r}+i+\ell_{(r-1)}%
-\ell_{(\beta-1)}<i_{\beta-1}.$ Once again, the inequalities for the lower
indices in the expression of $b$ in first fundamental are strict inequalities
as they were for $a$ but whose number now is increased by one. Obviously we
will have that
\[
i_{\beta}=i_{r}+i+\ell_{(r-1)}-\ell_{(\beta-1)}=\epsilon_{k}%
\]
for certain $1\leq k\leq q.$ Apply the $k^{th}$ left-transfer to obtain
\[
a^{\prime}=\tau_{\ell}^{k}(T_{n+2})=d_{(\tilde{\imath},\tilde{\ell})}%
^{\tilde{k}}(b^{\prime}),
\]
where $\tilde{k}=\alpha-1,$ $\tilde{\imath}=i_{k}^{\prime}-i_{\alpha}^{\prime
},$\ $\tilde{\ell}=\ell_{k}^{\prime}$\ and $\alpha$ is the smallest integer
$k\leq\alpha\leq q+1$ with $i_{\alpha}^{\prime}\leq i_{k}^{\prime}.$ Then
$b^{\prime}\otimes b$ is a component of $\Delta T_{n+2}$ as well.

\underline{\textit{Case IIb.}} Let $i>0$ and $i_{r}+i+\ell_{(r-1)}%
-\ell_{(\beta-1)}=i_{\beta-1}.$ Apply the $(\beta-1)^{th}$ left-transfer to
obtain
\[
b=\tau_{\ell}^{\beta-1}(T_{n+2})=d_{(\tilde{\imath},\tilde{\ell})}^{\beta
-1}(c),
\]
where $\tilde{\imath}=0$ and $\tilde{\ell}={\ell}_{\beta-1}.$ Then $a^{\prime
}\otimes c$ is also a component of $\Delta T_{n+2}$.

\underline{\textit{Case IIc.}} Let $i=0;$ there are two subcases:

\underline{\textit{Subcase i.}} If $1\leq r\leq p$ and no integer $1\leq k\leq
q$ exists with $r=t_{k},$ then%
\[
\epsilon_{k}=i_{r}+\ell+\ell_{(r-1)}-\ell_{(o\left(  k\right)  )}%
\ \ \text{and}\ \ i_{k}^{\prime}=i_{r}-\ell_{(o^{\prime}(r))}^{\prime}.
\]
Apply the $(r+1)^{th}$ left-transfer to obtain
\[
b=\tau_{\ell}^{r+1}(T_{n+2})=d_{(\tilde{\imath},\tilde{\ell})}^{\tilde{r}%
}(c),
\]
where $\tilde{r}=\alpha,\ \tilde{\imath}=i_{r}-i_{\alpha}+\ell_{(\alpha
-1)}-\ell_{(r)},$ $\tilde{\ell}=\ell_{r}-\ell$ and\ $\alpha$ is the smallest
integer $r<\alpha\leq p+1$ such that $i_{\alpha}+\ell_{(\alpha)}-\ell
_{(r)}\geq i_{r};$ namely, $\alpha=r+1$ or $\alpha=t_{o^{\prime}(r)}.$ Now the
required inequality for the $i_{k}^{\prime}$'s could conceivably be violated
for $k>o^{\prime}(r),$ but this is not so since it is easy to see that: (a) if
$r$ is not realized as $t_{k}$ for some $k>o^{\prime}(r),$ then each $z$ with
$\alpha<z<r$ is so; and (b) if $r=t_{k}$ for some $k>o^{\prime}(r),$ while
$\epsilon_{k}=i_{r}+\ell+\ell_{(r-1)}-\ell_{(o\left(  k\right)  )},$ then
$\alpha$ (and not $r$) serves as $t_{k}$ for indices of face operators in
expressions of $a$ and $c$; moreover, for either $\alpha=r+1$ or
$\alpha=t_{o^{\prime}(r)}$ (in which case $i_{k}^{\prime}\leq i_{o^{\prime
}(r)}^{\prime}$), one has $i_{k}^{\prime}\leq\min\limits_{o^{\prime}\left(
\alpha\right)  <j<k}\left\{  i_{j}^{\prime},\text{ }i_{\alpha}-\ell
_{(o^{\prime}\left(  \alpha\right)  )}^{\prime}\right\}  $ so that $a^{\prime
}\otimes c$ is also a component of $\Delta T_{n+2}$.

\underline{\textit{Subcase ii.}} If $1\leq r\leq p+1$ and $r=t_{k}$ for some
$1\leq k\leq q,$ then%
\[
\epsilon_{k}=i_{r}+\ell+\ell_{(r-1)}-\ell_{(o\left(  k\right)  )}%
\ \ \text{and}\ \ i_{k}^{\prime}=i_{r}-\ell_{(o^{\prime}(r))}^{\prime},
\]
from which we establish the equality
\[
\ell+\ell_{(r-1)}-\ell_{(o(k))}=\ell_{(k)}^{\prime}-\ell_{(o(r))}^{\prime}.
\]
The hypotheses of Lemma \ref{l1} are satisfied by putting $k_{1}%
=r=t_{k},\ p_{{1}}=o(k)$ and $k_{2}=k,\ p_{{2}}=o^{\prime}(r)$ so that
$a^{\prime}\otimes c$ is a component of $\Delta d_{(\tilde{\imath},\tilde
{\ell})}(T_{n+2})$ with $\tilde{\imath}=i_{r}$ and $\tilde{\ell}=\ell
+\ell_{(r-1)}-\ell_{(o(k))}.$

\underline{\textit{Case III.}} Let $c^{\prime}\otimes c$ be a component of
$\Delta d_{(i,\ell)}(T_{n+2}).$ Reduce $c$ and $c^{\prime}$ to the first and
second fundamental forms, respectively. According to Lemma \ref{l1}, we have
either $i=i_{r_{1}}$ or $i=i_{r_{2}}^{\prime}+\ell_{r_{2}}^{\prime}%
+\cdots+\ell_{k_{2}}=i_{r_{2}+1}^{\prime}+\ell_{r_{2}+1}^{\prime}$ for certain
integers $r_{1},$ $r_{2},$ $k_{2},$ i.e., $i_{r_{1}+1}=i_{r_{1}}$ or
$i_{r_{2}+1}^{\prime}+\ell_{r_{2}+1}^{\prime}=i_{r_{2}}^{\prime}.$ Note that
the shuffles under consideration prevent both cases from occurring
simultaneously, and we obtain the situation dual to either Subcase ii of Case
IIc, or to Subcase ii of Case Ic. Thus we obtain components $c^{\prime}\otimes
a$ or $a^{\prime}\otimes c$ of $\Delta T_{n+2}$ with $d_{(i,\ell)}^{r_{1}%
}(a)=c$ or $d_{(i,\ell)}^{r_{2}}(a^{\prime})=c^{\prime},$ respectively.

\section{Application: Tensor Products of $A_{\infty}$-(co)algebras}

In this section, we use the diagonal $\Delta$ to define the tensor product of
$A_{\infty}$-(co)algebras in maximal generality.$\ $We note that a special
case was given by J. Smith \cite{Smith} for certain objects with a richer
structure than we have here.$\ $We also mention that Lada and Markl
\cite{Lada} defined an $A_{\infty}$ tensor product structure on a construct
different from the tensor product of graded modules.

We adopt the following notation and conventions: Let $R$ be a commutative ring
with unity; $R$-modules are assumed to be $\mathbb{Z}$-graded, tensor products
and $Hom$'s are defined over $R$ and all maps are $R$-module maps unless
otherwise indicated.$\ $If an $R$-module $V$ is connected, $\overline
{V}=V/V_{0}$.$\ $The symbol $1:V\rightarrow V$ denotes the identity map; the
suspension and desuspension maps $\uparrow$ and $\downarrow$ shift dimension
by $+1$ and $-1$, respectively.$\ $Define $V^{\otimes0}=R$ and $V^{\otimes
n}=V\otimes\cdots\otimes V$ with $n>0$ factors; then $TV=\oplus_{n\geq
0}V^{\otimes n}$ and $T^{a}V$ (respectively, $T^{c}V$) denotes the free tensor
algebra (respectively, cofree tensor coalgebra) of $V.\ $Given $R$-modules
$V_{1},\ldots,V_{n},$ a permutation $\sigma\in S_{n}$ induces an isomorphism
$\sigma:V_{1}\otimes\cdots\otimes V_{n}\rightarrow V_{\sigma^{-1}(1)}%
\otimes\cdots\otimes V_{\sigma^{-1}(n)}$ by $\sigma(x_{1}\cdots x_{n})=\pm$
$x_{\sigma^{-1}(1)}\cdots x_{\sigma^{-1}(n)},$ where $\pm$ is the Koszul sign.
In particular, $\sigma_{2,n}=\left(  1\text{ }3\text{ }\cdots\text{ }\left(
2n-1\right)  \text{ }2\text{ }4\text{ }\cdots\text{ }2n\right)  :\left(
A\otimes B\right)  ^{\otimes n}\rightarrow A^{\otimes n}\otimes B^{\otimes n}$
and $\sigma_{n,2}=\sigma_{2,n}^{-1}$ induce isomorphisms $\left(  \sigma
_{2,n}\right)  ^{\ast}:Hom\left(  A^{\otimes n}\otimes B^{\otimes n},A\otimes
B\right)  \rightarrow Hom\left(  \left(  A\otimes B\right)  ^{\otimes
n}\hspace*{-0.02in},A\otimes B\right)  $ and $\left(  \sigma_{n,2}\right)
_{\ast}$:$\,Hom\left(  A\otimes B,A^{\otimes n}\otimes B^{\otimes n}\right)
\rightarrow Hom\left(  A\otimes B,\right.  $ $\left.  \left(  A\otimes
B\right)  ^{\otimes n}\right)  $. The map $\iota:Hom(U,V)\otimes Hom\left(
U^{\prime},V^{\prime}\right)  \rightarrow Hom\left(  U\otimes U^{\prime
},V\otimes V^{\prime}\right)  $ is the canonical isomorphism. If $f:V^{\otimes
p}\rightarrow V^{\otimes q}$ is a map, we let $f_{i,n-p-i}=1^{\otimes
i}\otimes f\otimes1^{\otimes n-p-i}:V^{\otimes n}\rightarrow V^{\otimes
n-p+q},$ where $0\leq i\leq n-p.\ $The abbreviations \textit{DGM, DGA,}and
\textit{DGC} stand for \textit{differential graded }$R$\textit{-module, DG
}$R$-\textit{algebra }and \textit{DG }$R$-\textit{coalgebra}, respectively.

We begin with a review of $A_{\infty}$-(co)algebras paying particular
attention to the signs.$\ $Let $A$ be a connected $R$-module equipped with
operations $\{\varphi^{k}\in$\linebreak$Hom^{k-2}\left(  A^{\otimes
k},A\right)  \}_{k\geq1}.\ $For each $k$ and $n\geq1,$ linearly extend
$\varphi^{k}$ to $A^{\otimes n}$ via
\[
\sum\limits_{i=0}^{n-k}\varphi_{i,n-k-i}^{k}:A^{\otimes n}\rightarrow
A^{\otimes n-k+1},
\]
and consider the induced map of degree $-1$ given by
\[
\sum\limits_{i=0}^{n-k}\left(  \uparrow\varphi^{k}\downarrow^{\otimes
k}\right)  _{i,n-k-i}:\left(  \uparrow\overline{A}\right)  ^{\otimes
n}\rightarrow\left(  \uparrow\overline{A}\right)  ^{\otimes n-k+1}.
\]
Let $\widetilde{B}A=T^{c}\left(  \uparrow\overline{A}\right)  $ and define a
map $d_{\widetilde{B}A}:\widetilde{B}A\rightarrow\widetilde{B}A$ of degree
$-1$ by
\begin{equation}
d_{\widetilde{B}A}=\sum\limits_{\substack{1\leq k\leq n \\0\leq i\leq n-k
}}\left(  \uparrow\varphi^{k}\downarrow^{\otimes k}\right)  _{i,n-k-i}.
\label{nine}%
\end{equation}
The identities $\left(  -1\right)  ^{\left[  n/2\right]  }\uparrow^{\otimes
n}\downarrow^{\otimes n}=1^{\otimes n}$ and $\left[  n/2\right]  +\left[
\left(  n+k\right)  /2\right]  \equiv nk+\left[  k/2\right]  $ (mod 2) imply
that
\begin{equation}
d_{\widetilde{B}A}=\sum\limits_{\substack{1\leq k\leq n \\0\leq i\leq n-k
}}\left(  -1\right)  ^{\left[  \left(  n-k\right)  /2\right]  +i\left(
k+1\right)  }\uparrow^{\otimes n-k+1}\varphi_{i,n-k-i}^{k}\downarrow^{\otimes
n}. \label{d-tilde-bar}%
\end{equation}

\begin{definition}
$\left(  A,\varphi^{n}\right)  _{n\geq1}$ is an \underline{$A_{\infty}%
$-algebra} if $d_{\widetilde{B}A}^{2}=0.$
\end{definition}

\begin{proposition}
\label{relations}For each $n\geq1,$ the operations $\left\{  \varphi
^{n}\right\}  $ on an $A_{\infty}$-algebra satisfy the following quadratic
relations:
\begin{equation}
\sum_{\substack{0\leq\ell\leq n-1 \\0\leq i\leq n-\ell-1}}\left(  -1\right)
^{\ell\left(  i+1\right)  }\varphi^{n-\ell}\varphi_{i,n-\ell-1-i}^{\ell+1}=0.
\label{A-infty-alg}%
\end{equation}

\end{proposition}

\begin{proof}
For $n\geq1,$
\begin{align*}
0  &  =\sum\limits_{\substack{1\leq k\leq n \\0\leq i\leq n-k}}\left(
-1\right)  ^{\left[  \left(  n-k\right)  /2\right]  +i\left(  k+1\right)
}\uparrow\varphi^{n-k+1}\downarrow^{\otimes n-k+1}\uparrow^{\otimes
n-k+1}\varphi_{i,n-k-i}^{k}\downarrow^{\otimes n}\\
&  =\sum\limits_{\substack{1\leq k\leq n \\0\leq i\leq n-k}}\left(  -1\right)
^{n-k+i\left(  k+1\right)  }\varphi^{n-k-1}\varphi_{i,n-k-i}^{k}\\
&  =-\left(  -1\right)  ^{n}\sum\limits_{\substack{0\leq\ell\leq n-1 \\0\leq
i\leq n-\ell-1}}\left(  -1\right)  ^{\ell\left(  i+1\right)  }\varphi^{n-\ell
}\varphi_{i,n-\ell-1-i}^{\ell+1}.
\end{align*}

\end{proof}

\noindent It is easy to prove that

\begin{proposition}
If $\left(  A,\varphi^{n}\right)  _{n\geq1}$ is an $A_{\infty}$-algebra, then
$\left(  \widetilde{B}A,d_{\widetilde{B}A}\right)  $ is a DGC.
\end{proposition}

\begin{definition}
Let $\left(  A,\varphi^{n}\right)  _{n\geq1}$ be an $A_{\infty}$%
-algebra.$\ $\underline{The tilde bar construction on $A$} is the DGC $\left(
\widetilde{B}A,d_{\widetilde{B}A}\right)  .$
\end{definition}

\begin{definition}
Let $A$ and $C$ be $A_{\infty}$-algebras.$\ $A chain map $f=f^{1}:A\rightarrow
C$ is a \underline{\textit{map of }$A_{\infty}$\textit{-algebras}}%
\textit{\ }if there exists a sequence of maps $\{f^{k}\in Hom^{k-1}\left(
A^{\otimes k},C\right)  \}_{k\geq2}$ such that
\[
\widetilde{f}=\sum_{n\geq1}\left(  \sum_{k\geq1}\uparrow f^{k}\downarrow
^{\otimes k}\right)  ^{\otimes n}:\widetilde{B}A\rightarrow\widetilde{B}C
\]
is a DGC map.
\end{definition}

Dually, consider a sequence of operations $\{\psi^{k}\in Hom^{k-2}\left(
A,A^{\otimes k}\right)  \}_{k\geq1}.\ $For each $k$ and $n\geq1,$ linearly
extend each $\psi^{k}$ to $A^{\otimes n}$ via
\[
\sum\limits_{i=0}^{n-1}\psi_{i,n-1-i}^{k}:A^{\otimes n}\rightarrow A^{\otimes
n+k-1},
\]
and consider the induced map of degree $-1$ given by
\[
\sum\limits_{i=0}^{n-1}\left(  \downarrow^{\otimes k}\psi^{k}\uparrow\right)
_{i,n-1-i}:\left(  \downarrow\overline{A}\right)  ^{\otimes n}\rightarrow
\left(  \downarrow\overline{A}\right)  ^{\otimes n+k-1}.
\]
Let $\widetilde{\Omega}A=T^{a}\left(  \downarrow\overline{A}\right)  $ and
define a map $d_{\widetilde{\Omega}A}:\widetilde{\Omega}A\rightarrow
\widetilde{\Omega}A$ of degree $-1$ by
\[
d_{\widetilde{\Omega}A}=\sum\limits_{\substack{n,k\geq1 \\0\leq i\leq n-1
}}\left(  \downarrow^{\otimes k}\psi^{k}\uparrow\right)  _{i,n-1-i},
\]
which can be rewritten as
\begin{equation}
d_{\widetilde{\Omega}A}=\sum\limits_{\substack{n,k\geq1 \\0\leq i\leq n-1
}}\left(  -1\right)  ^{\left[  n/2\right]  +i\left(  k+1\right)  +k\left(
n+1\right)  }\downarrow^{\otimes n+k-1}\psi_{i,n-1-i}^{k}\uparrow^{\otimes n}.
\label{d-tilde-cobar}%
\end{equation}

\begin{definition}
$\left(  A,\psi^{n}\right)  _{n\geq1}$ is an \underline{$A_{\infty}%
$-coalgebra} if $d_{\widetilde{\Omega}A}^{2}=0.$
\end{definition}

\begin{proposition}
For each $n\geq1,$ the operations $\left\{  \psi^{k}\right\}  $ on an
$A_{\infty}$-coalgebra satisfy the following quadratic relations:
\begin{equation}
\sum_{\substack{0\leq\ell\leq n-1 \\0\leq i\leq n-\ell-1}} \left(  -1\right)
^{\ell\left(  n+i+1\right)  }\psi_{i,n-\ell-1-i}^{\ell+1}\psi^{n-\ell}=0.
\label{A-infty-coalg}%
\end{equation}

\end{proposition}

\begin{proof}
The proof is similar to the proof of Proposition \ref{relations} and is omitted.
\end{proof}

\noindent Again, it is easy to prove that

\begin{proposition}
If $\left(  A,\psi^{n}\right)  _{n\geq1}$ is an $A_{\infty}$-coalgebra , then
$\left(  \widetilde{\Omega}A,d_{\widetilde{\Omega}A}\right)  $ is a DGA.
\end{proposition}

\begin{definition}
Let $\left(  A,\psi^{n}\right)  _{n\geq1}$ be an $A_{\infty}$-coalgebra.$\ $%
The \underline{tilde cobar construction on} \underline{$A$}is the DGA $\left(
\widetilde{\Omega}A,d_{\widetilde{\Omega}A}\right)  .$
\end{definition}

\begin{definition}
Let $A$ and $B$ be $A_{\infty}$-coalgebras.$\ $A chain map $g=g^{1}%
:A\rightarrow B$ is a \underline{\textit{map of }$A_{\infty}$%
\textit{-coalgebras}}\textit{\ }if there exists a sequence of maps $\{g^{k}\in
Hom^{k-1}\left(  A,B^{\otimes k}\right)  \}_{k\geq2}$ such that
\begin{equation}
\widetilde{g}=\sum_{n\geq1}\left(  \sum_{k\geq1}\downarrow^{\otimes k}%
g^{k}\uparrow\right)  ^{\otimes n}:\widetilde{\Omega}A\rightarrow
\widetilde{\Omega}B, \label{dga}%
\end{equation}
is a DGA map.
\end{definition}

The structure of an $A_{\infty}$-(co)algebra is encoded by the quadratic
relations among its operations (also called \textquotedblleft higher
homotopies\textquotedblright). Although the \textquotedblleft
direction,\textquotedblright\ i.e., sign, of these higher homotopies is
arbitrary, each choice of directions determines a set of signs in the
quadratic relations, the \textquotedblleft simplest\textquotedblright\ of
which appears on the algebra side when no changes of direction are made; see
(\ref{nine}) and (\ref{A-infty-alg}) above.$\ $Interestingly, the
\textquotedblleft simplest\textquotedblright\ set of signs appear on the
coalgebra side when $\psi^{n}$ is replaced by $\left(  -1\right)  ^{\left[
\left(  n-1\right)  /2\right]  }\psi^{n},$ $n\geq1,$ i.e., the direction of
every third and fourth homotopy is reversed.$\ $The choices one makes will
depend on the application; for us the appropriate choices are as in
(\ref{A-infty-alg}) and (\ref{A-infty-coalg}).

Let $\mathcal{A}_{\infty}=\oplus_{n\geq2}C_{\ast}\left(  K_{n}\right)  $ and
let $\left(  A,\varphi^{n}\right)  _{n\geq1}$ be an $A_{\infty}$-algebra with
quadratic relations as in (\ref{A-infty-alg}).$\ $For each $n\geq2,$ associate
$e^{n-2}\in C_{n-2}\left(  K_{n}\right)  $ with the operation $\varphi^{n}$
via
\begin{equation}
e^{n-2}\mapsto\left(  -1\right)  ^{n}\varphi^{n} \label{two}%
\end{equation}
and each codimension $1$ face $d_{(i,\ell)}\left(  e^{n-2}\right)  \in
C_{n-3}\left(  K_{n}\right)  $ with the quadratic composition
\begin{equation}
d_{(i,\ell)}\left(  e^{n-2}\right)  \mapsto\varphi^{n-\ell}\varphi
_{i,n-\ell-1-i}^{\ell+1}. \label{three}%
\end{equation}
Then (\ref{two}) and (\ref{three}) induce a chain map
\begin{equation}
\zeta_{A}:\mathcal{A}_{\infty}\longrightarrow\oplus_{n\geq2}Hom^{\ast}\left(
A^{\otimes n},A\right)  \label{alg}%
\end{equation}
representing the $A_{\infty}$-algebra structure on $A$.$\ $Dually, if $\left(
A,\psi^{n}\right)  _{n\geq1}$ is an $A_{\infty}$-coalge- bra with quadratic
relations as in (\ref{A-infty-coalg}), the associations
\[
e^{n-2}\mapsto\psi^{n}\text{$\ $and$\ $}d_{(i,\ell)}\left(  e^{n-2}\right)
\mapsto\psi_{i,n-\ell-1-i}^{\ell+1}\psi^{n-\ell}%
\]
induce a chain map
\begin{equation}
\xi_{A}:\mathcal{A}_{\infty}\longrightarrow\oplus_{n\geq2}Hom^{\ast}\left(
A,A^{\otimes n}\right)  \label{coalg}%
\end{equation}
representing the $A_{\infty}$-coalgebra structure on $A.\ $The definition of
the tensor product is now immediate:

\begin{definition}
\label{tensor product}The \underline{tensor product of $A_{\infty}$-algebras}
$\left(  A,\zeta_{A}\right)  $ and $\left(  B,\zeta_{B}\right)  $ is given by%
\[
\left(  A,\zeta_{A}\right)  \otimes\left(  B,\zeta_{B}\right)  =\left(
A\otimes B,\zeta_{A\otimes B}\right)  ,
\]
where $\zeta_{A\otimes B}$ is the sum of the compositions
\[%
\begin{array}
[c]{ccc}%
C_{\ast}(K_{n}) & \overset{{\LARGE \zeta}_{A\otimes B}}{\longrightarrow} &
Hom((A\otimes B)^{\otimes n},A\otimes B)\\
&  & \\
_{\strut\strut{\LARGE \Delta}_{K}}\text{$\ $}\downarrow\text{$\ $\ } &  &
\text{$\ \ \ \ \ \ $}\uparrow\text{$\ (\sigma_{2,n})^{\ast}\iota$}\\
&  & \\
C_{\ast}(K_{n})\otimes C_{\ast}(K_{n}) & \underset{{\LARGE \zeta}%
_{A}{\LARGE \otimes\zeta}_{B}}{\longrightarrow} & Hom(A^{\otimes n},A)\otimes
Hom\left(  B^{\otimes n},B\right)
\end{array}
\]
over all $n\geq2;$ the $A_{\infty}$-algebra operations $\Phi^{n}$ on $A\otimes
B$ are given by%
\[
\Phi^{n}=\left(  \sigma_{2,n}\right)  ^{\ast}\text{$\iota$}\left(  \zeta
_{A}\otimes\zeta_{B}\right)  \Delta_{K}\left(  e^{n-2}\right)  .
\]
Dually, the \underline{tensor product of $A_{\infty}$-coalgebras} $\left(
A,\xi_{A}\right)  $ and $\left(  B,\xi_{B}\right)  $ is given by%
\[
\left(  A,\xi_{A}\right)  \otimes\left(  B,\xi_{B}\right)  =\left(  A\otimes
B,\xi_{A\otimes B}\right)  ,
\]
where $\xi_{A\otimes B}$ is the sum of the compositions
\[%
\begin{array}
[c]{ccc}%
C_{\ast}(K_{n}) & \overset{{\LARGE \xi}_{A\otimes B}}{\longrightarrow} &
Hom(A\otimes B,(A\otimes B)^{\otimes n})\\
&  & \\
_{\strut\strut{\LARGE \Delta}_{K}}\text{$\ $}\downarrow\text{$\ $\ } &  &
\text{$\ \ \ \ \ \ $}\uparrow\text{$\ (\sigma_{n,2})_{\ast}\iota$}\\
&  & \\
C_{\ast}(K_{n})\otimes C_{\ast}(K_{n}) & \underset{{\LARGE \xi}_{A}%
{\LARGE \otimes\xi}_{B}}{\longrightarrow} & Hom(A,A^{\otimes n})\otimes
Hom\left(  B,B^{\otimes n}\right)
\end{array}
\]
over all $n\geq2;$ the $A_{\infty}$-coalgebra operations $\Psi^{n}$ on
$A\otimes B$ are given by
\[
\Psi^{n}=\left(  \sigma_{n,2}\right)  _{\ast}\text{$\iota$}\left(  \xi
_{A}\otimes\xi_{B}\right)  \Delta_{K}\left(  e^{n-2}\right)  .
\]

\end{definition}

\begin{example}
\label{iterate}If $\left(  A,\psi^{n}\right)  _{n\geq1}$ is an $A_{\infty}%
$-coalgebra,$\ $the $A_{\infty}$ operations on $A\otimes A$ are:%
\[%
\begin{array}
[c]{ll}%
\Psi^{1}= & \psi^{1}\otimes1+1\otimes\psi^{1}\\
& \\
\Psi^{2}= & \sigma_{2,2}\left(  \psi^{2}\otimes\psi^{2}\right)  \\
& \\
\Psi^{3}= & \sigma_{3,2}\left(  \psi_{0}^{2}\psi_{0}^{2}\otimes\psi^{3}%
+\psi^{3}\otimes\psi_{1}^{2}\psi_{0}^{2}\right)  \\
& \\
\Psi^{4}= & \sigma_{4,2}\left(  \psi_{0}^{2}\psi_{0}^{2}\psi_{0}^{2}%
\otimes\psi^{4}+\psi^{4}\otimes\psi_{2}^{2}\psi_{1}^{2}\psi_{0}^{2}+\psi
_{0}^{3}\psi_{0}^{2}\otimes\psi_{1}^{2}\psi_{0}^{3}\right.  \\
& \hspace{0.3in}%
\begin{array}
[c]{l}%
\left.  +\psi_{0}^{3}\psi_{0}^{2}\otimes\psi_{1}^{3}\psi_{0}^{2}+\psi_{1}%
^{2}\psi_{0}^{3}\otimes\psi_{1}^{3}\psi_{0}^{2}-\psi_{0}^{2}\psi_{0}%
^{3}\otimes\psi_{2}^{2}\psi_{0}^{3}\right)
\end{array}
\\
& \\
\Psi^{5}= & \sigma_{5,2}\left(  \psi_{0}^{2}\psi_{0}^{2}\psi_{0}^{2}\psi
_{0}^{2}\otimes\psi^{5}+\psi^{5}\otimes\psi_{3}^{2}\psi_{2}^{2}\psi_{1}%
^{2}\psi_{0}^{2}\right.  \\
& \hspace{0.3in}%
\begin{array}
[c]{l}%
+\psi_{0}^{4}\psi_{0}^{2}\otimes\psi_{2}^{3}\psi_{1}^{2}\psi_{0}^{2}-\psi
_{1}^{2}\psi_{0}^{4}\otimes\psi_{3}^{2}\psi_{1}^{3}\psi_{0}^{2}+\psi_{0}%
^{3}\psi_{0}^{2}\psi_{0}^{2}\otimes\psi_{1}^{3}\psi_{0}^{2}\\
+\psi_{0}^{4}\psi_{0}^{2}\otimes\psi_{2}^{2}\psi_{1}^{3}\psi_{0}^{2}-\psi
_{1}^{3}\psi_{0}^{3}\otimes\psi_{2}^{2}\psi_{1}^{3}\psi_{0}^{2}+\psi_{0}%
^{3}\psi_{0}^{2}\psi_{0}^{2}\otimes\psi_{1}^{4}\psi_{0}^{2}\\
+\psi_{0}^{4}\psi_{0}^{2}\otimes\psi_{2}^{2}\psi_{1}^{2}\psi_{0}^{3}-\psi
_{1}^{3}\psi_{0}^{3}\otimes\psi_{2}^{3}\psi_{1}^{2}\psi_{0}^{2}-\psi_{0}%
^{2}\psi_{0}^{3}\psi_{0}^{2}\otimes\psi_{2}^{2}\psi_{0}^{4}\\
+\psi_{0}^{2}\psi_{0}^{4}\otimes\psi_{3}^{2}\psi_{2}^{2}\psi_{0}^{3}+\psi
_{0}^{2}\psi_{0}^{2}\psi_{0}^{3}\otimes\psi_{3}^{2}\psi_{0}^{4}+\psi_{0}%
^{2}\psi_{0}^{3}\psi_{0}^{2}\otimes\psi_{2}^{3}\psi_{0}^{3}\\
+\psi_{2}^{2}\psi_{0}^{4}\otimes\psi_{2}^{3}\psi_{1}^{2}\psi_{0}^{2}+\psi
_{0}^{2}\psi_{1}^{2}\psi_{0}^{3}\otimes\psi_{2}^{3}\psi_{0}^{3}+\psi_{1}%
^{2}\psi_{0}^{3}\psi_{0}^{2}\otimes\psi_{1}^{3}\psi_{0}^{3}\\
-\psi_{0}^{3}\psi_{0}^{3}\otimes\psi_{3}^{2}\psi_{1}^{2}\psi_{0}^{3}+\psi
_{1}^{2}\psi_{1}^{2}\psi_{0}^{3}\otimes\psi_{1}^{4}\psi_{0}^{2}+\psi_{1}%
^{2}\psi_{0}^{3}\psi_{0}^{2}\otimes\psi_{1}^{4}\psi_{0}^{2}\\
-\psi_{0}^{3}\psi_{0}^{3}\otimes\psi_{3}^{2}\psi_{1}^{3}\psi_{0}^{2}+\psi
_{0}^{3}\psi_{0}^{2}\psi_{0}^{2}\otimes\psi_{1}^{2}\psi_{0}^{4}),
\end{array}
\end{array}
\]
\hspace*{0.4in}etc.
\end{example}

Note that the compositions in Definition \ref{tensor product} only use the
operations $\psi^{n}$ and not the quadratic relations (\ref{A-infty-coalg}).
Indeed, one can iterate an arbitrary family of operations $\left\{  \psi
^{n}\right\}  $ as in Example (\ref{iterate}) to produce iterated operations
$\Psi^{n}:A^{\otimes k}\rightarrow\left(  A^{\otimes k}\right)  ^{\otimes n}$
whether or not $\left(  A,\psi^{n}\right)  $ is an $A_{\infty}$-coalgebra. Of
course, the $\Psi^{n}$'s define an $A_{\infty}$-coalgebra structure on
$A^{\otimes k}$ whenever $d_{\widetilde{\Omega}\left(  A^{\otimes k}\right)
}^{2}=0,$ and we make extensive use of this fact in the sequel \cite{SU3}.
Finally, since $\Delta_{K}$ is homotopy coassociative (not strict), the tensor
product only iterates up to homotopy. In the sequel we always coassociate on
the extreme left.

\section{Appendix: Associahedral Sets}

An associahedral set is a combinatorial object generated by Stasheff
associahedra $K$ and equipped with appropriate face and degeneracy operators.
Associahedral sets are similar in many ways to simplicial or cubical sets.

\subsection{Singular associahedral sets}

To motivate the notion of an associahedral set, we begin with a construction
of singular associahedral sets, our universal example. Let $X$ be a
topological space. Define the singular associahedral complex $Sing^{K}X$ as
follows: Let
\[
({Sing}^{K}X)_{n-k+2}^{(j_{1},n_{1}),...,(j_{k+1},n_{k+1})}=\{\text{Continuous
maps }K_{n_{1}+2}\times\cdots\times K_{n_{k+1}+2}\rightarrow X\},
\]
where $\sum_{q=1}^{k+1}n_{q}=n-k,$\thinspace\ $n_{q}\geq0,\,0\leq k\leq
n,\,n\geq j_{1}\geq\cdots\geq j_{k}\geq j_{k+1}=0$ and $K_{n_{1}+2}%
\times\cdots\times K_{n_{k+1}+2}$ is a Cartesian product of associahedra. Let
\[%
\begin{tabular}
[c]{l}%
${\delta^{q}}_{(i_{q},\ell_{q})}:K_{n_{1}+2}\times\cdots\times K_{n_{\beta
-1}+2}\times K_{\ell_{q}+1}\times\cdots$\\
$\hspace*{1in}\hspace*{1in}\cdots\times K_{n_{q}-\ell_{q}+2}\times
K_{n_{q+1}+2}\times\cdots\times K_{n_{k+2}+2}$\\
$\hspace*{1in}\rightarrow K_{n_{1}+2}\times\cdots\times K_{n_{q}+2}%
\times\cdots K_{n_{k+2}+2},$%
\end{tabular}
\ \ \ \ \ \
\]
be the map determined by $\delta_{(i_{q},\ell_{q})}^{q}={\delta^{\prime}%
}_{(i_{q},\ell_{q})}^{q}\circ T_{(q,\beta)},$ where
\[%
\begin{tabular}
[c]{l}%
${\delta^{\prime}}_{(i_{q},\ell_{q})}^{q}:K_{n_{1}+2}\times\cdots\times
K_{\ell_{q}+1}\times K_{n_{q}-\ell_{q}+2}\times\cdots\times K_{n_{k+2}+2}$\\
$\hspace*{1in}\rightarrow K_{n_{1}+2}\times\cdots\times K_{n_{q}+2}%
\times\cdots K_{n_{k+2}+2},$%
\end{tabular}
\ \
\]
${\delta^{\prime}}_{(i_{q},\ell_{q})}^{q}=1^{\times q-1}\times\delta
_{(i_{q},\ell_{q})}^{\prime\prime}\times1^{\times k+1-q},$
\[
{\delta^{\prime\prime}}_{(i_{q},\ell_{q})}^{q}:K_{\ell_{q}+1}\times
K_{n_{q}-\ell_{q}+2}\rightarrow K_{n_{q}+2}%
\]
is the standard inclusion corresponding to the pair $(i_{q},\ell_{q}),$ and
\[%
\begin{tabular}
[c]{l}%
$T_{(q,\beta)}:K_{n_{1}+2}\times\cdots\times K_{n_{\beta-1}+2}\times
K_{\ell_{q}+1}\times\cdots$\\
$\hspace*{1in}\hspace*{1in}\cdots\times K_{n_{q}-\ell_{q}+2}\times
K_{n_{q+1}+2}\times\cdots\times K_{n_{k+2}+2}$\\
$\hspace*{1in}\overset{\approx}{\rightarrow}K_{n_{1}+2}\times\cdots\times
K_{n_{\beta-1}+2}\times\cdots$\\
$\hspace*{1.5in}\cdots\times K_{\ell_{q}+1}\times K_{n_{q}-\ell_{q}+2}\times
K_{n_{q+1}+2}\times\cdots\times K_{n_{k+2}+2}$%
\end{tabular}
\ \
\]
is the permutation isomorphism in which $\beta$ is defined by (\ref{ksign2}).
Let
\[
\eta_{i}^{q}:K_{n_{1}+2}\times\cdots\times K_{n_{q}+3}\times\cdots
K_{n_{k+1}+2}\rightarrow K_{n_{1}+2}\times\cdots\times K_{n_{q}+2}\times
\cdots\times K_{n_{k+1}+2}%
\]
be the projection (cf. \cite{Stasheff}). Then for $f\in(Sing^{K}%
X)_{n-k+2}^{(j_{1},n_{1}),...,(j_{k+1},n_{k+1})},$ define%

\[%
\begin{array}
[c]{l}%
d_{\left(  i_{q},\ell_{q}\right)  }^{q}:(Sing^{K}X)_{n-k+2}^{(j_{1}%
,n_{1}),...,(j_{k+1},n_{k+1})}\rightarrow\\
\ \ \ \ \ \ \ \ \ \ \ \ \ \ \ \ \ \ \ \ \ (Sing^{K}X)_{n-k+1}^{(j_{1}%
,n_{1}),...,(j_{\beta}-1,n_{\beta}-1),({j(q,\beta)},\ell_{q}-1),...,(j_{q}%
,n_{q}-\ell_{q}),...,(j_{k+1},n_{k+1})},
\end{array}
\]
with $j(q,\beta)$ is defined in (\ref{ksign1}) and
\[%
\begin{array}
[c]{l}%
s_{i}^{q}:(Sing^{K}X)_{n-k+2}^{(j_{1},n_{1}),\ldots,(j_{k+1},n_{k+1}%
)}\rightarrow(Sing^{K}X)_{n-k+3}^{(j_{1},n_{1}),\ldots,(j_{q},n_{q}%
+1),\ldots,(j_{k+1},n_{k+1})},
\end{array}
\]
as compositions
\[
d_{(i_{q},\ell_{q})}^{q}(f)=f\circ\delta_{(i_{q},\ell_{q})}^{q}\text{ and
}s_{i}^{q}(f)=\eta_{i}^{q}\circ f.
\]
Given the abstract definition below, one can easily check that $(Sing^{K}%
X,\ d_{(i_{q},\ell_{q})},\ s_{i}^{q})$ is an associahedral set.

\subsection{Abstract associahedral sets}

\begin{definition}
\label{associahedral}An \underline{associahedral set} is a graded set
$\vspace*{0.1in}$\newline$\mathcal{K}=\{K_{n-k+2}^{(j_{1},n_{1}),...,(j_{k+1}%
,n_{k+1})}|\,n\geq j_{1}\geq\cdots\geq j_{k+1}=0,\,n_{q}\geq0,n_{\left(
k+1\right)  }=n-k\}_{0\leq k\leq n,}\vspace*{0.1in}$
\end{definition}

\noindent\textit{together with face and degeneracy operators defined for
}$1\leq q\leq k+1$:$\vspace*{0.1in}$

\noindent$d_{\left(  i_{q},\ell_{q}\right)  }^{q}:K_{n-k+2}^{(j_{1}%
,n_{1}),...,(j_{k+1},n_{k+1})}\rightarrow\vspace*{0.1in}$

\noindent$K_{n-k+1}^{(j_{1},n_{1}),...,(j_{\beta-1},n_{\beta-1}),({j(q,\beta
)},\ell_{q}-1),(j_{\beta},n_{\beta}),...,(j_{q-1},n_{q-1}),(j_{q},n_{q}%
-\ell_{q}),(j_{q+1},n_{q+1}),...,(j_{k+1},n_{k+1})}\vspace*{0.1in}$

\noindent\textit{where }$j(q,\beta)\ $\textit{and}$\ \beta\ $are
\textit{defined in (\ref{ksign1}) and (\ref{ksign2}), }$0\leq i_{q}\leq
n_{q};$ $1\leq\ell_{q}\leq n_{q};$ $i_{q}+\ell_{q}\leq n_{q}+1,$
and$\vspace*{0.1in}$

\noindent$s_{j}^{q}:K_{n-k+2}^{(j_{1},n_{1}),\ldots,(j_{k+1},n_{k+1}%
)}\rightarrow K_{n-k+3}^{(j_{1},n_{1}),\ldots,(j_{q-1},n_{q-1}),(j_{q}%
,n_{q}+1),(j_{q+1},n_{q+1}),\ldots,(j_{k+1},n_{k+1})}\vspace*{0.1in}$

\noindent\textit{for} $1\leq j\leq n_{q}+3,$ \textit{satisfying relations
(1)-(3) as well as}%
\[%
\begin{array}
[c]{ll}%
d_{(i,\ell)}^{p}s_{j}^{q}=s_{j}^{q+1}d_{(i,\ell)}^{p};\vspace{0.1in} & p<q\\
d_{(i,\ell)}^{p}s_{j}^{q}=s_{j}^{q}d_{(i,\ell)}^{p};\hspace*{0.4in}%
\vspace{0.1in} & p>q\\
d_{(i,\ell)}^{q}s_{j}^{q}=s_{j-\ell}^{q+1}d_{(i,\ell)}^{q};\vspace{0.1in} &
i+\ell+1<j\\
d_{(i,\ell)}^{q}s_{j}^{q}=s_{j-i}^{q}d_{(i,\ell-1)}^{q};\vspace{0.1in} &
i<j<i+\ell+2,\text{ }\ell>1\\
d_{(i,\ell)}^{q}s_{j}^{q}=s_{j}^{q+1}d_{(i-1,\ell)}^{q};\vspace{0.1in} & i\geq
j,\text{ }\ell\leq n_{q}\\
d_{(i,\ell)}^{q}s_{j}^{q}=1;\vspace{0.1in} & (i,\ell)=(j-1,1),\text{ }1\leq
j<n_{q}+3\\
d_{(i,\ell)}^{q}s_{j}^{q}=1;\vspace{0.1in} & (i,\ell)=(j-2,1),\text{ }1<j\leq
n_{q}+3\\
d_{(i,\ell)}^{q}s_{j}^{q}=1;\vspace{0.1in} & (i,\ell)=(0,n_{q}+1),\text{
}j=n_{q}+3\\
d_{(i,\ell)}^{q}s_{j}^{q}=1; & (i,\ell)=(1,n_{q}+1),\text{ }j=1\\
s_{j}^{p}s_{j^{\prime}}^{q}=s_{j^{\prime}}^{q}s_{j}^{p};\vspace{0.1in} & p\neq
q\\
s_{j}^{q}s_{j^{\prime}}^{q}=s_{j^{\prime}+1}^{q}s_{j}^{q};\vspace{0.1in} &
p=q,\text{ }j\leq j^{\prime}.
\end{array}
\]

Given an associahedral set $\mathcal{K},$ let%
\[
(C_{\ast}(\mathcal{K}),d)=\bigoplus C_{n-k}(K_{n-k+2}^{(j_{1},n_{1}%
),...,(j_{k+1},n_{k+1})},d^{n_{1},...,n_{k+1}}),
\]
where
\[
d^{n_{1},...,n_{k+1}}=\sum_{(i_{q},\ell_{q})}(-1)^{\epsilon_{1}+\epsilon_{2}%
}d_{(i_{q},\ell_{q})}^{q},
\]
with $\epsilon_{i}$ defined in (\ref{esign}); define the diagonal
\[
\Delta_{\mathcal{K}}:C_{\ast}(\mathcal{K})\rightarrow C_{\ast}(\mathcal{K}%
)\otimes C_{\ast}(\mathcal{K})
\]
on $K^{(0,n)}$ by (\ref{shuffle2}) and extend to other components of
$K^{(j_{1};n_{1}),...,(j_{k+1};n_{k+1})}$ by the formal multiplicative rule
with respect to indices $(n_{1},...,n_{k+1}),$ i.e., by the same formulas as
on a product cell $K_{n_{1}+2}\times\cdots\times K_{n_{k+1}+2}$. Finally, set
$C_{\ast}^{N}(\mathcal{K})=C_{\ast}(\mathcal{K})/D,$ where $D$ is the
submodule generated by the degeneracies; then $\left(  C_{\ast}^{N}%
(\mathcal{K}),d\right)  $ is a chain complex equipped with a diagonal
$\Delta_{\mathcal{K}}$ induced by $\Delta$.

\end{document}